\PassOptionsToPackage{numbers,sort&compress}{natbib}
\documentclass[smallextended,natbib]{svjour3}
\bibpunct[, ]{[}{]}{,}{n}{,}{,}

\usepackage[utf8x]{inputenc}
\usepackage[T1]{fontenc}
\usepackage{lmodern}
\usepackage[hidelinks,hypertexnames=false]{hyperref}
\usepackage[margin=1in]{geometry}
\usepackage{multicol}
\usepackage{booktabs}
\usepackage{amsfonts} \usepackage{amssymb} 
\usepackage{graphicx} % support the \includegraphics command and options
\usepackage{epstopdf}
 \usepackage{enumerate}
\usepackage{url} \usepackage{verbatim} \usepackage{float} 

\usepackage{tabularx,ragged2e,booktabs,caption} \usepackage{makecell}
\usepackage{algorithm2e}
\usepackage{enumitem} 
\usepackage{caption}

\usepackage{float}

\usepackage{makecell}
\usepackage{mathrsfs}
\usepackage{tikz}
\usepackage{multirow}
\usepackage{chngcntr}
\usepackage{mathtools}

\usepackage{subcaption}
\counterwithin{table}{subsection}

\DeclareMathOperator*{\relation}{Supp}
\DeclareMathOperator*{\partition}{Part}

\DeclareMathOperator*{\rank}{rank}

\DeclareMathOperator*{\trace}{Tr}
\DeclareMathOperator*{\Span}{span}

\DeclareMathOperator*{\nnz}{nnz}

\DeclareMathOperator*{\vect}{vect}

\newcommand{\norm}[1]{\lVert #1 \rVert }

\newcommand{\ip}[2]{#1\cdot#2}

\newcommand{\dualaff}{C+\mathcal{L}^{\perp}}
\newcommand{\primaff}{Y+\mathcal{L}}

\newcommand{\charmats}{\mathcal{B}_{\mathcal{P}}}

\newcommand{\coneNameTwo}{\mathcal{C}}
\newcommand{\algName}{\mathbb{S}^n}
\newcommand{\coneName}{\mathbb{S}^n_{+}}
\newcommand{\isoName}{\Psi}
\newcommand{\isoname}{\isoName}

\newcommand{\ones}{11^T}

\title{Dimension  reduction for semidefinite programs via Jordan algebras}

\author{Frank Permenter \and Pablo A. Parrilo} % <-this % stops a space
\begin{document}
\newtheorem{thm}{Theorem}
\numberwithin{thm}{section}

\newtheorem{cor}{Corollary}
\numberwithin{cor}{section}

\newtheorem{lem}{Lemma}
\numberwithin{lem}{section}

\newtheorem{rem}{Remark}

\newtheorem{prop}{Proposition}
\numberwithin{prop}{section}

\newtheorem{defn}{Definition}
\numberwithin{defn}{section}

\newtheorem{ex}{Example}
\numberwithin{ex}{section}

\newtheorem{ass}{Condition}
\numberwithin{ass}{section}

\newtheorem{prob}{Open Problem}

\date{\today}
\maketitle

%\setcounter{tocdepth}{4}
%\tableofcontents

\begin{abstract}
We propose a new method for simplifying semidefinite programs (SDP) inspired by
symmetry reduction. Specifically, we show if an orthogonal projection map
satisfies certain invariance conditions, restricting to its range yields an
equivalent primal-dual pair over a lower-dimensional symmetric cone---namely,
the cone-of-squares of a Jordan subalgebra of symmetric matrices. We present
a simple algorithm for minimizing the rank of this projection and hence the
dimension of this subalgebra. We also show that minimizing rank optimizes the
direct-sum decomposition of the algebra into simple ideals, yielding an
optimal ``block-diagonalization'' of the SDP.\@ Finally, we give
combinatorial versions of our algorithm that execute at reduced computational
cost  and illustrate effectiveness of an implementation on examples.  Through
the theory of Jordan algebras, the proposed method easily extends to linear
and second-order-cone programming and, more generally, symmetric cone
optimization. 
\end{abstract}

\section{Introduction}

Many practically relevant optimization problems can be posed as  semidefinite
programs (SDPs)---convex optimization problems over the cone of  positive
semidefinite (psd) matrices.  While SDPs are efficiently solved in theory,
specific instances may be intractable in practice unless one exploits special
structure.  Existing techniques for structure exploitation include   facial
reduction~\cite{borwein1981regularizing,drusvyatskiy2017many,pataki2013simple}
and symmetry reduction
\cite{symSoS_Gatermann200495,vallentin2009symmetry,de2010exploiting,bachoc2012invariant}.
In this paper, we present a method that builds on this latter technique.

To explain, we first recall a key step in symmetry reduction: finding  an
orthogonal projection map whose range intersects the solution set.  This
projection (called a \emph{Reynolds} or \emph{group-average} operator) maps
feasible points to feasible points without changing the objective function,
which implies that its range (called the \emph{fixed-point subspace}) contains
solutions.   This leads to a simple statement of our method: minimize
rank---or, equivalently, the dimension of the range---over a tractable subset
of maps with this property.  As we show, this minimization problem is
efficiently solved for arbitrary SDP instances by a simple algorithm. Further,
the subset of projections considered strictly contains those implicit in
existing symmetry reduction procedures (Section~\ref{sec:suffgen}); hence, our
method is more general.

Symmetry reduction not only reduces the dimension of the feasible set, it also
simplifies the semidefinite constraint. This simplification process is
informally called block-diagonalization, and it amounts to  finding a canonical
direct-sum decomposition of the fixed-point subspace.  The projection we
identify enables similar simplifications.  Precisely, the range is always a
subalgebra of the \emph{Jordan algebra} of real, symmetric matrices and hence
also has a canonical direct-sum decomposition into \emph{simple ideals}.  Further, its
intersection with the psd cone (the \emph{cone-of-squares} of the subalgebra)
has a corresponding decomposition into irreducible  \emph{symmetric
cones}~\cite[Chapter 3]{faraut1994analysis}. As we review (Section~\ref{sec:starDecomp}), finding
this decomposition generalizes current block-diagonalization
techniques based on *-algebras~\cite{de2011numerical,maehara2010numerical}. As we show, 
minimizing the rank of the projection  optimizes this decomposition in a
precise sense. 

Finally, our method easily extends to any symmetric cone optimization problem
(including linear and second-order-cone programs).  Indeed, via Jordan algebra
theory, our algorithm for finding projections extends ``word-by-word'',
mirroring similar  extensions of interior-point
methods~\cite{alizadeh2000symmetric}.  

We organize this paper as follows.    Section~\ref{sec:equivsdp} contains
preliminaries.  Section~\ref{sec:findingsub} gives an  algorithm for finding a
minimum-rank projection.    Section~\ref{sec:optdirectsum} shows that
minimizing rank yields an algebra with an optimal direct-sum decomposition. Section~\ref{sec:comb} gives
combinatorial (but less powerful) versions of our algorithm that can be less
costly to execute.  Computational results appear in Section~\ref{sec:example}.

\section{Preliminaries}\label{sec:equivsdp}

We  consider  a primal-dual pair of semidefinite programs (SDPs)   expressed in conic form (\cite[Chapter 4]{nesterov1994interior}):
\begin{align}\label{sdp:main}
\begin{array}{ll}
\mbox{ minimize } & \ip{C}{X} \\
\mbox{ subject to } & X \in Y + \mathcal{L}    \\
& X \in \mathbb{S}_{+}^n  \\
\end{array} \qquad
\begin{array}{ll}
\mbox{ maximize } & -\ip{Y}{S}    \\
\mbox{ subject to } & S \in C + \mathcal{L}^{\perp}   \\
& S \in \mathbb{S}_{+}^n.   \\
\end{array}
\end{align}
Here, $X \in \mathbb{S}^n$ and $S \in \mathbb{S}^n$ are  decision variables  in
the vector space  $\mathbb{S}^n$  of real symmetric matrices equipped with
trace inner-product $\ip{X}{Y}:= \trace XY$,  $\mathbb{S}_{+}^n \subseteq \mathbb{S}^n$
denotes the (self-dual) cone of psd matrices,     $\mathcal{L}
\subseteq \mathbb{S}^n$  is a linear subspace with orthogonal complement
$\mathcal{L}^{\perp}\subseteq \mathbb{S}^n$, and  $Y + \mathcal{L}$  and $C +
\mathcal{L}^{\perp}$  are affine sets defined by fixed $C \in \mathbb{S}^n$ and $Y\in \mathbb{S}^n$. 
We refer to $X$ and $S$ as the primal and dual decision variables, respectively, 
noting that we have identified the dual space ${(\mathbb{S}^n)}^*$ with $\mathbb{S}^n$.
(Note that in this form, the complementary slackness condition $\ip{X}{S} = 0$ does not
necessarily imply the primal and dual objective values are equal. Rather, they differ by a constant
that depends on the particular choice of $C$ and $Y$.)

Throughout this paper we also, for a subspace  $\mathcal{S} \subseteq \mathbb{S}^n$, let
$P_\mathcal{S}: \mathbb{S}^n \rightarrow \mathbb{S}^n$ denote the
corresponding orthogonal projection map, i.e., the unique self-adjoint
and idempotent map with range equal to $\mathcal{S}$.

\subsection{Constraint Set Invariance}
Our goal is to find a    subspace  $\mathcal{S} \subseteq \mathbb{S}^n$ that
contains primal and dual solutions of~\eqref{sdp:main} if they exist. To do
this, we will find a projection that maps
feasible points to feasible points without changing the cost function (which
implies the range contains solutions),  a key idea from  symmetry reduction~\cite{symSoS_Gatermann200495,vallentin2009symmetry,de2010exploiting,bachoc2012invariant}.
Precisely, we will search over the orthogonal projections that satisfy the following set of conditions,
which we'll show are also implicit in existing symmetry reduction approaches (Section~\ref{sec:suffgen}). 

\begin{defn}[Constraint Set Invariance Conditions]\label{cond:CI}
We say the orthogonal projection map $P_{\mathcal{S}} :
  \mathbb{S}^n \rightarrow \mathbb{S}^n$ satisfies the \emph{Constraint
Set Invariance Conditions} for the primal-dual pair~\eqref{sdp:main} if
  \begin{enumerate}[label= (\alph*),series=lafter] 
    \item\label{itm:posmap} $P_\mathcal{S}(\mathbb{S}^n_{+}) \subseteq  \mathbb{S}^n_{+}$, 
    i.e., $P_\mathcal{S}$ is a \emph{positive} map;
    \item\label{itm:primaff} $P_\mathcal{S} (Y+\mathcal{L}) \subseteq Y + \mathcal{L}$;
    \item\label{itm:dualaff}  $P_\mathcal{S} (C+ \mathcal{L}^{\perp}) \subseteq C + \mathcal{L}^{\perp}$.
  \end{enumerate}
\end{defn}
Under these conditions,  $P_\mathcal{S} :
\mathbb{S}^n \rightarrow \mathbb{S}^n$ maps primal/dual feasible points
to primal/dual feasible points (by definition).  For $C$ and all $X \in Y+\mathcal{L}$, these conditions also imply that 
\begin{align}\label{eq:orthCSI}
X-P_{\mathcal{S}}(X) \in \mathcal{L}  \qquad C-P_{\mathcal{S}}(C) \in \mathcal{L}^{\perp},
\end{align}
which in turn implies that $P_\mathcal{S}$ preserves the cost function on the primal feasible set:
\[
C \cdot X = P_{\mathcal{S}}(C) \cdot P_{\mathcal{S}}(X) = C \cdot P_S P_S(X) = C \cdot  P_{\mathcal{S}}(X).
\]
(Here, the first equality holds given~\eqref{eq:orthCSI}, and the second and  third given that $P_{\mathcal{S}}$
is  self-adjoint and idempotent.)  A similar argument   shows $Y \cdot S=Y \cdot  P_{\mathcal{S}}(S)$ for all dual feasible $S$. In summary, we've proven the following. 
\begin{prop}[Preservation of optimal values]\label{prop:optval}
        Suppose $P_{\mathcal{S}} : \mathbb{S}^n \rightarrow \mathbb{S}^n$
        satisfies the \emph{Constraint Set Invariance Conditions} for the
        primal-dual pair~\eqref{sdp:main}.  Let $\theta_p := \inf_{} \left\{ \ip{C}{X} :   X \in  \mathbb{S}^n_{+} \cap (Y  +
        \mathcal{L})  \right\}$ and $\theta_d :=   \sup_{}
        \left\{ -\ip{Y}{S} :   S \in  \mathbb{S}^n_{+} \cap (C  +
        \mathcal{L}^{\perp}) \right\}$. Then, 
                \begin{align*}
		\theta_p =\inf_{ } \left\{ \ip{C}{X} :   X \in  \mathbb{S}^n_{+} \cap (Y  + \mathcal{L}) \cap \mathcal{S} \right\}, \qquad
		\theta_d = \sup_{ } \left\{ -\ip{Y}{S} :   S \in  \mathbb{S}^n_{+} \cap (C  + \mathcal{L}^{\perp}) \cap \mathcal{S} \right\}.
		\end{align*}
                Further, $\mathcal{S}$ contains points that attain $\theta_{p}$ and $\theta_{d}$ when such points exist.
\end{prop}
In other words, we've shown that restricting the primal/dual feasible set to $\mathcal{S}$ doesn't change the primal/dual optimal value or its attainment.
\subsubsection{Infeasibility certificates}
A dual improving direction is  $S \in \mathbb{S}^n_{+} \cap
\mathcal{L}^{\perp}$ satisfying $Y \cdot S < 0$. Analogously, a primal
improving direction is  $X \in \mathbb{S}^n_{+} \cap \mathcal{L}$
satisfying $C \cdot X < 0$.  The existence of primal (resp.\ dual) improving directions implies
infeasibility of the dual (resp.\ primal).  It turns out that if
$P_{\mathcal{S}} : \mathbb{S}^n \rightarrow \mathbb{S}^n$ satisfies the
Constraint Set Invariance Conditions, then the subspace $\mathcal{S}$ contains
improving directions whenever they exist for the original problem.  To show this, we need the following lemma.
 \begin{lem}[Invariance of linear subspaces]\label{lem:invariantSub}
        Suppose $P_{\mathcal{S}} : \mathbb{S}^n \rightarrow \mathbb{S}^n$
        satisfies the Constraint Set Invariance Conditions.  Then $\mathcal{L}$
        and $\mathcal{L}^{\perp}$ are both invariant subspaces of
        $P_{\mathcal{S}}$, i.e.,   $P_{\mathcal{S}}(\mathcal{L}) \subseteq
        \mathcal{L}$ and $P_{\mathcal{S}}(\mathcal{L}^{\perp}) \subseteq
        \mathcal{L}^{\perp}$.
 	
        \begin{proof}
          For all $Z \in \mathcal{L}$, we have that $P_{\mathcal{S}}(Z) =
          P_{\mathcal{S}}(Y) - P_{\mathcal{S}}(Y-Z) \in \mathcal{L}$, where 
          containment in $\mathcal{L}$ follows given that $Y+\mathcal{L}$
          contains both  $P_{\mathcal{S}}(Y)$ and $P_{\mathcal{S}}(Y-Z)$  by
          the Constraint Set Invariance Conditions.  This shows that $\mathcal{L}$
          is an invariant subspace; the proof for $\mathcal{L}^{\perp}$ is identical.
        \end{proof}
 \end{lem}
We can now show the desired result.
\begin{prop}[Improving directions]\label{prop:ray}
Suppose $P_{\mathcal{S}} : \mathbb{S}^n \rightarrow \mathbb{S}^n$ satisfies the
Constraint Set Invariance Conditions.  The following statements hold.
	\begin{itemize}
          \item If $S \in \mathbb{S}^n$ is a dual improving direction, then so is $P_{\mathcal{S}}(S)$.
          \item If $X \in \mathbb{S}^n$ is a primal improving direction, then so is $P_{\mathcal{S}}(X)$.
\end{itemize}
  \begin{proof}
    Let $S$ be a dual improving direction.  Lemma~\ref{lem:invariantSub}  and the
    Constraint Set Invariance Conditions imply that $\mathbb{S}^n_{+} \cap
    \mathcal{L}^{\perp}$ contains $P_{\mathcal{S}}(S)$, that
    $\mathcal{L}^{\perp}$ contains $S-P_{\mathcal{S}}(S)$ and that
    $\mathcal{L}$ contains $Y - P_{\mathcal{S}}(Y)$.  These latter two facts
    imply that $S \cdot Y =   S \cdot P_{\mathcal{S}}(Y)$; hence,
    $P_{\mathcal{S}}(S)$ is a dual improving direction.  Proof of the second
    statement is identical.
  \end{proof}
\end{prop}

\subsubsection{Restricted primal-dual pair}

We've seen that intersecting the primal and dual feasible with $\mathcal{S}$
does not change the primal and dual optimal value if $P_{\mathcal{S}}$
satisfies the Constraint Set Invariance Conditions (Proposition~\ref{prop:optval}). Further,   $\mathcal{S}$
contains solutions or infeasibility certificates for~\eqref{sdp:main} when
such objects exists (Propositions~\ref{prop:optval}-\ref{prop:ray}). These facts allow us to solve~\eqref{sdp:main} by first
restricting the primal and dual to $\mathcal{S}$.  The following shows that these
restrictions are a primal-dual pair if we view $\mathcal{S}$ as the ambient
space. 
\begin{prop}[Duality and restrictions]\label{prop:invar}
	Suppose that $P_{\mathcal{S}} : \mathbb{S}^n \rightarrow \mathbb{S}^n$ satisfies  the Constraint Set Invariance Conditions (Definition~\ref{cond:CI}). Then, treating the range $\mathcal{S}$ as the ambient space, the pair  of optimization problems 
	\begin{align}\label{sdp:equiv}
	\begin{array}{ll}
	\mbox{ \rm minimize } & \ip{P_\mathcal{S}(C)}{X} \\
	\mbox{ \rm subject to } & X \in P_\mathcal{S}(Y) + \mathcal{L} \cap \mathcal{S}     \\
	& X \in \mathbb{S}_{+}^n \cap \mathcal{S}  \\
	\end{array} \qquad
	\begin{array}{ll}
	\mbox{ \rm maximize } & -\ip{P_\mathcal{S}(Y)}{S} \\
	\mbox{ \rm subject to } & S \in P_\mathcal{S}(C) +  \mathcal{L}^{\perp} \cap \mathcal{S}   \\
	& S \in \mathbb{S}_{+}^n \cap \mathcal{S} \\
	\end{array}
	\end{align}
	is a primal-dual pair, i.e., 
	\begin{align}\label{eq:ambient}
          {(\mathbb{S}_{+}^n \cap \mathcal{S})}^{*} \cap \mathcal{S} &= \mathbb{S}_{+}^n \cap \mathcal{S},   \qquad	{(\mathcal{L} \cap \mathcal{S})}^{\perp} \cap   \mathcal{S} = \mathcal{L}^{\perp} \cap \mathcal{S}.
	\end{align} 
	Moreover,
	\begin{align}\label{eq:feasEq}
	(Y+\mathcal{L}) \cap \mathcal{S} =	P_\mathcal{S}(Y) + \mathcal{L} \cap \mathcal{S}, \qquad  (C + \mathcal{L}^{\perp}) \cap \mathcal{S} = P_\mathcal{S}(C) + \mathcal{L}^{\perp} \cap \mathcal{S}.
	\end{align}
	\newcommand{\projName}{P_{\mathcal{S}}}
	\begin{proof}
		
		For any set $\mathcal{T} \subseteq \mathbb{S}^n$, the condition $P_{\mathcal{S}}(\mathcal{T}) \subseteq \mathcal{T}$ 	 implies $P_{\mathcal{S}}(\mathcal{T}) = \mathcal{S} \cap \mathcal{T}$ given that $P_{\mathcal{S}}$ is the orthogonal projection onto $\mathcal{S}$. Using this fact, we have that
		\[
		(Y + \mathcal{L}) \cap \mathcal{S} = 	\projName(Y + \mathcal{L}) = \projName(Y) + \projName(\mathcal{L}) =  \projName(Y) +  \mathcal{L} \cap \mathcal{S},
		\]
                where the last equality uses the additional fact that
                $P_\mathcal{S}(\mathcal{L}) \subseteq \mathcal{L}$
                (Lemma~\ref{lem:invariantSub}).  The other equality
                in~\eqref{eq:feasEq} follows by identical argument. The
                inclusions $\supseteq$ of~\eqref{eq:ambient} are obvious. To see
                the inclusions $\subseteq$,  let $\mathcal{T}$ be any set
                satisfying $P_{\mathcal{S}}(\mathcal{T}) \subseteq
                \mathcal{T}$. Then, for any $X \in {(\mathcal{T} \cap
                \mathcal{S})}^* \cap \mathcal{S}$,
		\[
		\langle  X,  Y \rangle  = \langle P_\mathcal{S}(X),  Y \rangle = \langle   X, P_\mathcal{S}(Y) \rangle \ge 0, \qquad 	\forall Y \in \mathcal{T},
		\]
		where the first equality holds since $X \in \mathcal{S}$, the second equality since $P_\mathcal{S}$ is self-adjoint and the inequality since $P_\mathcal{S}(Y) \in \mathcal{T} \cap \mathcal{S}$. Hence, $X \in \mathcal{T}^*$.
	\end{proof}

\end{prop}
We illustrate this proposition with the following example.
\begin{ex}
	Consider the following primal-dual  pair of semidefinite programs:
	\begin{align*}
	\begin{array}{lll}
	\mbox{\rm minimize  } & x_1+x_2      & \\
	\mbox{\rm subject to} \\
	&  \left(\begin{array}{ccccc} 
	x_1 & 1 & x_3 &  x_4 \\ 
	1 & x_2 & x_4 & -x_3 \\
	x_3 & x_4  &1 &  x_5 \\ 
	x_4 & -x_3 & x_5 & 0  
	\end{array}\right)  \succeq 0 
	\end{array} 
	\begin{array}{lll}
	\mbox{\rm maximize } & -(s_5+2s_{1}) & \\
	\mbox{\rm {subject to}} \\
	& \left(\begin{array}{ccccc} 
	1 & s_{1} & s_{2} & s_{3} \\ 
	s_{1} & 1 & -s_{3}& s_{2}\\
	s_{2} & -s_{3} & s_5  & 0 \\ 
	s_{3} & s_{2} & 0  &   s_6 
	\end{array}\right) \succeq 0.
	\end{array}
	\end{align*}
        The projection   $P_{\mathcal{S}} : \mathbb{S}^4 \rightarrow \mathbb{S}^4$ satisfies the Constraint Set Invariance Conditions~(Definition~\ref{cond:CI}) if $\mathcal{S}$ equals the span of $\{E_{21}+E_{12} \} \cup {\{ E_{ii} \}}^3_{i=1}$.   Hence,  
	one obtains primal and dual optimal solutions by solving the  following restrictions to $\mathcal{S}$:
	\begin{align*}
	\begin{array}{lll}
	\mbox{\rm minimize  } & x_1+x_2     & \\
	\mbox{\rm subject to} \\
	& \left(\begin{array}{ccccc} 
	x_1 & 1 & 0 &  0 \\ 
	1 & x_2 & 0 & 0 \\
	0& 0  &1 &  0 \\ 
	0 & 0 & 0 &0  
	\end{array}\right) \succeq 0
	\end{array} 
	\begin{array}{lll}
	\mbox{\rm maximize } & -(s_5+2s_{1})  & \\
	\mbox{\rm {subject to}} \\
	& \left(\begin{array}{ccccc} 
	1 & s_{1} & 0 & 0 \\ 
	s_{1} & 1 & 0& 0\\
	0 & 0 & s_5  & 0 \\ 
	0 & 0 & 0  &   0
	\end{array}\right) \succeq 0.
	\end{array}
	\end{align*}	
\end{ex}

\subsubsection{Relationship with prior work}\label{sec:suffgen}

A common symmetry reduction technique, described in~\cite{symSoS_Gatermann200495}, assumes existence of a subgroup  $\mathcal{G} \subset \mathbb{R}^{n \times n}$ of
the group of orthogonal matrices 
that, for all $U \in \mathcal{G}$, satisfies
\begin{align}\label{eq:GroupCondition} 
  UCU^T = C,   \qquad \left\{ U X U^T   : X \in Y+\mathcal{L} \right\} \subseteq Y+\mathcal{L}.   
\end{align} 
Under this condition, one can restrict the primal problem to the \emph{fixed-point subspace} 
\begin{align}\label{eq:commdefn} 
\mathcal{F}_{\mathcal{G}} := \{ X \in \mathbb{S}^{n}  : UXU^T = X \;\; \forall U \in \mathcal{G} \}, 
\end{align} without changing its optimal value~\cite[Theorem 3.3]{symSoS_Gatermann200495}, in analogy
with Proposition~\ref{prop:optval}. (One can also derive analogues of Proposition~\ref{prop:invar} based on these conditions; see, e.g.,~\cite[Proposition~2]{dobre2015exploiting}.)
It turns out that the orthogonal projection onto $\mathcal{F}_{\mathcal{G}}$ (called
the \emph{Reynolds operator}) satisfies the Constraint Set Invariance Conditions. 

To see this, first observe that  $P_{\mathcal{F}_\mathcal{G}}$ is the map $X \mapsto
\frac{1}{|\mathcal{G}|} \sum_{U \in \mathcal{G}} U X U^T$.
As shown in~\cite{symSoS_Gatermann200495},
\begin{align}\label{eq:reynoldprop}
  P_{\mathcal{F}_{\mathcal{G}}}(\mathbb{S}^n_{+}) \subseteq \mathbb{S}^n_{+}, \qquad P_{\mathcal{F}_{\mathcal{G}}}(Y + \mathcal{L}) \subseteq Y+\mathcal{L}, \qquad P_{\mathcal{F}_{\mathcal{G}}}(C)=C,
\end{align}
 given~\eqref{eq:GroupCondition} and the fact that $P_{\mathcal{F}_{\mathcal{G}}}(X)$ is a convex
combination of points in $\{ U X U^T : U \in \mathcal{G}\}$.  The proof of the next proposition shows that $P_{\mathcal{F}_{\mathcal{G}}}$ also satisfies $P_{\mathcal{F}_{\mathcal{G}}}(C+\mathcal{L}^{\perp}) \subseteq C+
\mathcal{L}^{\perp}$  (and hence the full set of Constraint Set Invariance
Conditions). 
\begin{prop}[Constraint Set Invariance From
	Groups]\label{prop:csigroup} Let $\mathcal{G} \subset \mathbb{R}^{n \times n}$ be a finite group of orthogonal matrices that satisfies~\eqref{eq:GroupCondition}.  Then, $P_{\mathcal{F}_{\mathcal{G}}} : \mathbb{S}^n \rightarrow \mathbb{S}^n$ satisfies the Constraint Set Invariance Conditions (Definition~\ref{cond:CI}), and, in
	addition, the equality $P_{\mathcal{F}_{\mathcal{G}}}(C) = C$.
	
        \begin{proof} Given~\eqref{eq:reynoldprop},  we only need to show that $P_{\mathcal{F}_{\mathcal{G}}}(C+\mathcal{L}^{\perp}) \subseteq \mathcal{L}^{\perp}$.  To begin, $P_{\mathcal{F}_{\mathcal{G}}}(Y + \mathcal{L})
		\subseteq Y + \mathcal{L}$ implies that $P_{\mathcal{F}_{\mathcal{G}}}(\mathcal{L}) \subseteq \mathcal{L}$ (Lemma~\ref{lem:invariantSub}). Since $P_{\mathcal{F}_{\mathcal{G}}}$ is self-adjoint, $P_{\mathcal{F}_{\mathcal{G}}}(\mathcal{L}) \subseteq \mathcal{L}$
		holds if and only if $P_{\mathcal{F}_{\mathcal{G}}}(\mathcal{L}^{\perp}) \subseteq  \mathcal{L}^{\perp}$. Since $P_{\mathcal{F}_{\mathcal{G}}}(C) = C$, we conclude that $P_{\mathcal{F}_{\mathcal{G}}}(C+\mathcal{L}^{\perp}) \subseteq
		C+\mathcal{L}^{\perp}$, as desired.
	\end{proof}
\end{prop}

Another technique, surveyed in~\cite{de2010exploiting},  treats $\mathbb{R}^{n \times n}$ as a \emph{*-algebra} with matrix multiplication as a product  and transposition as a *-involution.  It then finds any  \emph{*-subalgebra}, i.e., any subspace closed under matrix multiplication and transposition, that contains the primal affine set $Y + \mathcal{L}$. If $\mathcal{M} \subseteq \mathbb{R}^n$ is such a *-subalgebra, then $\mathcal{S} := \mathcal{M} \cap \mathbb{S}^n$ contains primal and dual solutions~\cite[Theorem 2]{de2010exploiting}.
Further, the projection $P_{\mathcal{S}}$ satisfies
\begin{align}\label{eq:starAlgprop}
P_{\mathcal{S}}(\mathbb{S}^n_{+}) \subseteq  \mathbb{S}^n_{+},  \qquad P_{\mathcal{S}}(Y + \mathcal{L}) = Y + \mathcal{L},
\end{align}
where the inclusion $P_{\mathcal{S}}(\mathbb{S}^n_{+}) \subseteq
\mathbb{S}^n_{+}$ holds because $\mathcal{M}$ is a *-subalgebra. It turns out
that $P_{\mathcal{S}_{}}(C+\mathcal{L}^{\perp}) \subseteq C+\mathcal{L}^{\perp}$ (and hence the full set of Constraint Set Invariance
Conditions) also holds.
\begin{prop}[Constraint Set Invariance From *-algebras]\label{prop:staralgdata}
	Let  $\mathcal{M} \subseteq \mathbb{R}^{n \times n}$ be any *-subalgebra containing   $Y + \mathcal{L} \subseteq \mathbb{S}^n$. Let $\mathcal{S} =   \mathcal{M} \cap \mathbb{S}^n$.  Then $P_{\mathcal{S}} : \mathbb{S}^n \rightarrow \mathbb{S}^n$ satisfies the Constraint Set Invariance Conditions (Definition~\ref{cond:CI}), and, in addition, the equality $P_{\mathcal{S}}(Y + \mathcal{L}) = Y + \mathcal{L}$.
	\begin{proof}
          Given~\eqref{eq:starAlgprop}, we only need to show that $P_{\mathcal{S}_{}}(C+\mathcal{L}^{\perp}) \subseteq C+\mathcal{L}^{\perp}$.

To begin, since $\mathcal{S}$ contains $Y+\mathcal{L}$, we have
that $P_{\mathcal{S}} (Y+\mathcal{L}) = Y+\mathcal{L}$, which in turn implies that
  \begin{align}\label{starAlg:proof1}
    \mathcal{L} &=  P_{\mathcal{S}} (\mathcal{L}).
  \end{align}
From \eqref{starAlg:proof1}, we conclude that $\mathcal{L}$ is an invariant
          subspace of $P_{\mathcal{S}}$ which in turn in implies that
          $\mathcal{L}^{\perp}$ is an invariant subspace of the adjoint of
          $P_{\mathcal{S}}$. But  $P_{\mathcal{S}}$ is self-adjoint; hence, 
  \begin{align}\label{starAlg:proof2}
    P_{\mathcal{S}} (\mathcal{L}^{\perp}) \subseteq \mathcal{L}^{\perp}.
  \end{align}
We conclude that
\[
  P_{\mathcal{S}} (C_\mathcal{\mathcal{L}} + \mathcal{L}^{\perp}) =
          C_\mathcal{\mathcal{L}} +P_{\mathcal{S}} ( \mathcal{L}^{\perp})   \subseteq
          C_\mathcal{\mathcal{L}} + \mathcal{L}^{\perp},
\]
where the equality follows from  \eqref{starAlg:proof1} and the inclusion from \eqref{starAlg:proof2}.  Since $C_{\mathcal{L}} + \mathcal{L}^{\perp} = C + \mathcal{L}^{\perp}$, the conclusion follows.
\end{proof}
\end{prop}
Note that this proposition puts no condition on objective matrix $C$ of the primal problem.
Similarly, \cite[Theorem 2]{de2010exploiting} puts no condition on the dual objective function.

\subsection{Reformulations over isomorphic, symmetric cones}

The fixed-point subspace of symmetry reduction and *-subalgebras  have
structured intersections with $\mathbb{S}^n_{+}$: each intersection is
isomorphic to a direct product of psd cones of Hermitian matrices with real, complex,
or quaternion entries.  Such a product is an instance of a \emph{symmetric
cone}, a special type of self-dual cone that admits efficient optimization
algorithms~\cite{faybusovich1997linear,alizadeh2000symmetric}.   To maintain this feature, Section~\ref{sec:unitalCond} gives an
additional condition on the projection $P_{\mathcal{S}} : \mathbb{S}^n \rightarrow \mathbb{S}^n$ that ensures 
$\mathbb{S}^n_{+} \cap \mathcal{S}$ is always isomorphic to a symmetric cone.  This next proposition shows that the primal-dual
pair~\eqref{sdp:main} can be explicitly reformulated over such an isomorphic cone under the Constraint
Set Invariance Conditions.

  \begin{prop}[Reformulations over isomorphic cones]\label{prop:equiv2}
          Suppose $P_{\mathcal{S}} : \mathbb{S}^n \rightarrow \mathbb{S}^n$
          satisfies the Constraint Set Invariance Conditions
          (Definition~\ref{cond:CI}).   For an inner-product space $\mathcal{V}$,
          let $\isoName : \mathcal{V} \rightarrow \mathbb{S}^n$ be an injective
          linear map with range equal to $\mathcal{S}$ and $\coneNameTwo
          \subseteq \mathcal{V}$ a self-dual cone  that satisfies
          \begin{align} \label{eq:coneEqBig}
          \isoname(\coneNameTwo) = \mathbb{S}^n_{+} \cap \mathcal{S}.
          \end{align}
If $\hat X \in \mathcal{V}$ 
          and $\hat S \in \mathcal{V}$ solve the  primal-dual pair of conic optimization problems
          \begin{align}\label{sdp:equiv2}
          \begin{array}{ll}
          \mbox{ \rm minimize } &  \langle \isoname^*(C), \hat X  \rangle \\
            \mbox{ \rm subject to } & \hat X \in {(\isoname^*\isoname)}^{-1}  \isoname^* (Y +\mathcal{L})    \\
          & \hat X \in \coneNameTwo  
          \end{array} \qquad
          \begin{array}{ll}
            \mbox{ \rm maximize } & -  \langle {(\isoname^*\isoname)}^{-1}  \isoname^*(Y), \hat S \rangle  \\
          \mbox{ \rm subject to } & \hat S \in \isoname^*( C +\mathcal{L}^{\perp} )    \\
          & \hat S \in \coneNameTwo,
          \end{array}
          \end{align}
        then  $\isoname(\hat X)$ and $\isoname {(\isoname^*\isoname)}^{-1}(\hat S)$ solve the primal-dual pair~\eqref{sdp:equiv}---and hence the primal-dual pair~\eqref{sdp:main}.
          
          \begin{proof}
            We will show that $\isoname$ and ${(\isoname^*\isoname)}^{-1}  \isoname^*$ are mappings between primal feasible points of~\eqref{sdp:equiv2} and~\eqref{sdp:equiv}
            that do not change the objective value. To see that $\isoname$ has this property, let $\hat X$ be a feasible point of~\eqref{sdp:equiv2}. Then,
            \[
              \isoname( \hat X ) \in P_{\mathcal{S}}( Y+\mathcal{L}), \qquad \isoname( \hat X ) \in \mathbb{S}^n_{+} \cap \mathcal{S}, \qquad  \langle C, \isoname(\hat X) \rangle = \langle \isoname^*(C), \hat X \rangle,
            \]
            where the first containment follows given that
            $\isoname{(\isoname^* \isoname)}^{-1} \isoname^*$ equals
            $P_{\mathcal{S}}$ and the second by~\eqref{eq:coneEqBig}. Since $P_{\mathcal{S}}(Y+\mathcal{L}) \subseteq Y+\mathcal{L}$, we conclude that $\isoname(\hat X)$ is feasible for~\eqref{sdp:equiv} with same objective as $\hat X$.

            Now suppose $X$ is feasible for~\eqref{sdp:equiv}. Then $X = \isoname(\hat X)$ for a unique
            $\hat X \in \coneNameTwo$ since $\isoname$ is injective.  Indeed,
            we must have that $\hat X = {(\isoname^*\isoname)}^{-1}  \isoname^*
            (X)$, since
            \[
              \isoname {(\isoname^*\isoname)}^{-1}  \isoname^* (X) = P_{\mathcal{S}} (X) = X.
            \]
            Hence, $\hat X = {(\isoname^*\isoname)}^{-1}  \isoname^* (X)$ is a feasible point of~\eqref{sdp:equiv2} with objective
            \[
            \langle \isoname^*(C), \hat X  \rangle = \langle \isoname^*(C), {(\isoname^*\isoname)}^{-1}  \isoname^* X \rangle = \langle C, \isoname  {(\isoname^*\isoname)}^{-1}  \isoname^* X \rangle = \langle C,  P_{\mathcal{S}} X\rangle = \langle C,   X\rangle,
            \]
            as desired. 
                 
            For the dual, we similarly prove that $\isoname^*$ and $\isoname {(\isoname^*\isoname)}^{-1}$ are mappings between the feasible
            sets that do not change the objective. The proof is almost the same, but exploits the additional fact that
            \begin{align}\label{eq:conesEq}
            \isoname^* \isoname(\coneNameTwo) = \coneNameTwo,
            \end{align}
            which we show in the Appendix (Lemma~\ref{lem:sdimage}). To begin,  if $\hat S$ is dual feasible for~\eqref{sdp:equiv2}, then $\isoname {(\isoname^*\isoname)}^{-1}( \hat S )$ 
            satisfies
                          \[
                            \isoname {(\isoname^*\isoname)}^{-1}( \hat S ) \in  P_{\mathcal{S}}(C+\mathcal{L}^{\perp}),   \quad 	\isoname {(\isoname^*\isoname)}^{-1}( \hat S ) \in \mathbb{S}^n_{+} \cap \mathcal{S}, \quad \langle Y, \isoname {(\isoname^*\isoname)}^{-1}( \hat S )\rangle =  \langle {(\isoname^*\isoname)}^{-1}  \isoname^*(Y), \hat S \rangle,
                          \]
            Here, the first containment follows because $\hat S \in \isoname^*(C+\mathcal{L}^{\perp})$
            and $\isoname{(\isoname^*\isoname)}^{-1} \isoName^* = P_{\mathcal{S}}$; the second  by~\eqref{eq:coneEqBig}~and~\eqref{eq:conesEq}.
            Since $P_{\mathcal{S}}(C+\mathcal{L}^{\perp}) \subseteq C+\mathcal{L}^{\perp}$, we conclude that $\isoname{(\isoname^*\isoname)}^{-1} (\hat S)$ is feasible for~\eqref{sdp:equiv} with same objective as $\hat S$.
            
            On the other hand, if $S$ is dual
            feasible for~\eqref{sdp:equiv}, then $\hat S :=   \isoname^* S$
            must be the unique $\hat S \in \mathcal{V}$ satisfying $S =
            \isoname {(\isoname^*\isoname)}^{-1}    \hat S$ since $\isoname
            {(\isoname^*\isoname)}^{-1} \isoname^* S = S$. Further,  $\hat S
            \in \coneNameTwo$ since $\isoName^* (\mathbb{S}^n_{+} \cap \mathcal{S}) = \coneNameTwo$ by \eqref{eq:coneEqBig}~and~\eqref{eq:conesEq}.
            Hence $\hat S$ is dual feasible for \eqref{sdp:equiv2}. Further, its objective satisfies
            \[
              \langle {(\isoname^*\isoname)}^{-1}\isoname^*(Y), \hat S \rangle = \langle {(\isoname^*\isoname)}^{-1}\isoname^*(Y), \isoname^*S \rangle =  \langle Y, P_{\mathcal{S}} S \rangle = \langle Y, S \rangle,
            \]
            as desired.
          \end{proof}
\end{prop}

\subsection{Euclidean Jordan Algebras}

We now develop a condition that guarantees $\mathbb{S}^n_{+} \cap \mathcal{S}$
is isomorphic to a symmetric cone $\coneNameTwo$ and discuss how to construct a
linear map $\isoname$ satisfying $\isoname(\coneNameTwo) = \mathbb{S}^n_{+}
\cap \mathcal{S}$. For this, we first view  $\mathbb{S}^n$ as a \emph{Euclidean
Jordan algebra}, i.e., as an inner-product space $\mathcal{J}$ equipped with
bilinear product $(x,y) \mapsto x\circ y$ that is commutative,  satisfies the
\emph{Jordan identity}
\[
  x^2 \circ(y \circ x) = (x^2 \circ y) \circ x  \qquad \forall x, y \in \mathcal{J}
\]
(where $x^2 := x \circ x$), and, for all fixed $x$, is a self-adjoint, i.e., $\langle x \circ y, z \rangle
= \langle  y, x \circ z \rangle$ for all $y,z \in \mathcal{J}$.  To satisfy
these axioms, we equip $\mathbb{S}^n$ with the trace
inner-product $X \cdot S := \trace XY$ and product $X\circ Y :=
\frac{1}{2}(XY+YX)$.  For any Euclidean Jordan algebra $\mathcal{J}$, the set of squares $\{ x^2 : x \in
\mathcal{J}\}$ is always a
symmetric cone~\cite[Chapter 3]{faraut1994analysis}, often called the
\emph{cone-of-squares} of $\mathcal{J}$.  For the aforementioned product, the cone-of-squares of
$\mathbb{S}^n$ is just the psd cone $\mathbb{S}^n_{+}$.

It turns out that $\mathbb{S}^n_{+} \cap \mathcal{S}$ is isomorphic to a
symmetric cone whenever $\mathcal{S}$  is a \emph{subalgebra} of
$\mathbb{S}^n$, i.e., whenever $\mathcal{S}$ contains $X \circ Y$ for all
$X,Y \in \mathcal{S}$.  This follows because $\mathcal{S}$ 
satisfies the Euclidean-Jordan-algebra axioms (when viewed as the ambient
space) and has cone-of-squares $\mathbb{S}^n_{+} \cap \mathcal{S}$.  As a
consequence, we can write $\mathbb{S}^n_{+} \cap \mathcal{S}$  as the linear image
$\isoName(\coneNameTwo)$ 
of the cone-of-squares $\coneNameTwo$ of any isomorphic algebra $\mathcal{J}$ using an injective \emph{homomorphism}
$\isoName : \mathcal{J} \rightarrow \mathbb{S}^n$, i.e., an injective linear map  satisfying $\isoname(x \circ y) = \isoname(x) \circ
\isoname(y)$. Formally: 
\begin{prop}\label{prop:coneOfSquareRep}
  Let $\mathcal{S}$ be a subalgebra of $\mathbb{S}^n$. Let $\mathcal{J}$ be any Euclidean Jordan algebra
  isomorphic to $\mathcal{S}$
  with cone-of-squares $\mathcal{C} \subseteq \mathcal{J}$. Let $\isoName : \mathcal{J} \rightarrow \mathbb{S}^n$ be an injective
  homomorphism with range equal to $\mathcal{S}$. Then, 
	\[
           \isoname(\coneNameTwo) = \mathbb{S}^n_{+} \cap \mathcal{S}.
	\]
\end{prop}
\begin{ex}
	Let $\mathcal{S}$ denote the subalgebra of $\mathbb{S}^n$ spanned by
	\[
	E_1 = \begin{bmatrix}
	1   &  0  &   0  &  0 \\
	0  &   1  &  0  & 0 \\
	0  &  0  &   0  &   0 \\
	0  & 0  &   0  &   0
	\end{bmatrix},  
	E_2 = \begin{bmatrix}
	0   &  0  &   0  &  0 \\
	0  &   0  &  0  & 0 \\
	0  &  0  &   1  &   0 \\
	0  & 0  &   0  &   1
	\end{bmatrix},  
	T_1 = \begin{bmatrix}
	0   &  0  &   1  &  0 \\
	0  &   0  &  0  & 1 \\
	1  &  0  &   0  &   0 \\
	0  & 1  &   0  &   0
	\end{bmatrix},  
	T_2 = \begin{bmatrix}
	0   &  0  &   0  &  1 \\
	0  &   0  &  -1  & 0 \\
	0  &  -1  &   0  &   0 \\
	1  & 0  &   0  &   0
	\end{bmatrix}.
	\]
      Let $\mathcal{J}$ denote the spin-factor algebra $\mathbb{R} \times \mathbb{R}^3$ with 
      cone-of-squares $\mathcal{Q}:=\{ (x_0, x) \in \mathcal{J} : \|x\|_2 \le x_0  \}$ and product
      \[
(x_0 , x) \circ  (y_0 , y) := (x_0 y_0 + x^{T}y, x_0y+y_0 x).
      \]
        Finally, let $\isoName : \mathcal{J}  \rightarrow \mathbb{S}^n$ denote the injective  linear map satisfying
	\[
	\isoName e_i  = E_i, \qquad \isoName t_i  =T_i \qquad i \in \{1,2\},
	\]
	where
	\[
	e_1 = \begin{bmatrix}
	\frac{1}{2}   &  \frac{1}{2}   &   0  &  0
	\end{bmatrix} ^T,  
	e_2 = \begin{bmatrix}
	\frac{1}{2}   &  -\frac{1}{2}   &   0  &  0
	\end{bmatrix} ^T,  
	t_1 = \begin{bmatrix}
	0   &  0  &   1  &  0 
	\end{bmatrix} ^T,  
	t_2 = \begin{bmatrix}
	0   &  0  &   0  &  1 
	\end{bmatrix} ^T.
	\]
        Then, the image of $\mathcal{Q}$ under $\isoName$ is  $\mathbb{S}^n_{+} \cap \mathcal{S}$.  Further, $\isoName$ is an injective homomorphism from $\mathcal{J}$ into $\mathbb{S}^n$ with range equal to $\mathcal{S}$.
\end{ex}
Given $\mathcal{S}$, one can find a canonical isomorphic algebra $\mathcal{J}$  and an injective homomorphism $\isoname$
numerically~\cite[Chapter 6]{permenterThesis2017}; see
Section~\ref{sec:subalgstruc} for more details.

\subsubsection{Positive projections, unitality, and subalgebras} \label{sec:unitalCond}

We can guarantee that $\mathcal{S}$ is a subalgebra and hence that $\mathbb{S}^n
\cap \mathcal{S}$ is isomorphic to a symmetric cone by revisiting the
positivity constraint $P_{\mathcal{S}}(\mathbb{S}^n_{+}) \subseteq
\mathbb{S}^n_{+}$ of the Constraint Set Invariance Conditions. Specifically, we
obtain this guarantee  by imposing positivity and, in addition, \emph{unitality}.
\begin{defn}[Unitality Condition]\label{defn:unital}
	We say that $P_{\mathcal{S}}: \mathbb{S}^n \rightarrow \mathbb{S}^n$ is unital if the range $\mathcal{S} $ contains a unit $E \in \mathcal{S}$ for the Jordan product $X \circ Y= \frac{1}{2}(XY+YX)$, i.e., if there exists $ E \in \mathcal{S}$ for which $X \circ E = X$ for all $X \in \mathcal{S}$.
\end{defn}
\begin{thm}[Characterization of positive, unital projections]\label{thm:posjor}
	Let $\mathcal{S} \subseteq \mathbb{S}^n$ be subspace with orthogonal projection map $P_{\mathcal{S}} : \mathbb{S}^n \rightarrow \mathbb{S}^n$.
	The following statements are equivalent.
	\begin{enumerate}
		\item The projection $P_{\mathcal{S}}$ is positive (Definition~\ref{cond:CI}-\ref{itm:posmap}) and unital (Definition~\ref{defn:unital}).
		\item The subspace $\mathcal{S}$ is subalgebra of $\mathbb{S}^n$.
	\end{enumerate}
	\begin{proof}
		See appendix.
	\end{proof}
\end{thm}
As we show in the appendix, this theorem follows from  basic linear algebra and arguments of \citet{stormer2012positive} (who proves an analogous result for complex Jordan algebras).
Note also that the unitality condition holds when $\mathcal{S}$ arises from a group or a *-subalgebra via Proposition~\ref{prop:csigroup}~or~\ref{prop:staralgdata}.
\subsubsection{Structure of subalgebras}\label{sec:subalgstruc}

We now discuss the structure of subalgebras in more detail. To begin, call an
abstract Euclidean Jordan algebra $\mathcal{J}$  \emph{simple} if its only
\emph{ideals} are $\{0\}$ and $\mathcal{J}$, where an ideal $\mathcal{I}
\subseteq \mathcal{J}$ is a subspace satisfying $x \circ y \in \mathcal{I}$ for
all $y \in \mathcal{J}$ and $x \in \mathcal{I}$.  Similarly, call an ideal
simple if it is simple when viewed as an algebra.  It is well known that any
subalgebra $\mathcal{S}$ equals an orthogonal direct-sum of its simple ideals
(e.g.,~\cite[Proposition~III.4.4]{faraut1994analysis}). Further, the
isomorphism classes of these ideals are fully understood~\cite[Chapter
V]{faraut1994analysis}. As a consequence,  $\mathcal{S}$ is always, up-to
linear transformation, a direct-sum of ``canonical'' algebras. Formally:
\begin{prop}[Structure theorem for subalgebras~\cite{faraut1994analysis}]\label{prop:structure}
  Let $\oplus^r_{i=1} \mathcal{S}_i$ be the orthogonal direct-sum decomposition of a subalgebra
  $\mathcal{S}\subseteq \mathbb{S}^n$ into simple ideals.  Then, there exists simple Jordan algebras $\mathcal{J}_1,\ldots,\mathcal{J}_r$
  and injective homomorphisms $\isoName_i : \mathcal{J}_i \rightarrow \mathbb{S}^n$ satisfying
  \begin{align}\label{eq:decompForm}
  \mathcal{S}_i = \isoName_i(\mathcal{J}_i),
  \end{align}
  where each $\mathcal{J}_i$ is one of the following\footnote{We omit
  the Albert algebra from this list since it is \emph{exceptional}, i.e., it is an algebra that is not \emph{special}. By definition, all subalgebras of $\mathbb{S}^n$
  are special \cite[2.3.1]{hanche1984jordan}; hence, no subalgebra of $\mathbb{S}^n$ is isomorphic to the Albert algebra.}
  \begin{itemize}
    \item A spin-factor algebra $\mathbb{R} \times \mathbb{R}^m$ with product $(x_0 , x) \circ  (y_0 , y) := (x_0 y_0 + x^{T}y, x_0y+y_0 x)$
    \item The set of Hermitian matrices $H_d(\mathcal{D})$ of order $d$ with entries in $\mathcal{D}$ and product $\frac{1}{2}(XY+YX)$, where $\mathcal{D}$
      denotes the real numbers, the complex numbers, or the quaternions.
  \end{itemize}
\end{prop}
Efficient algorithms exist for finding the ideals $\mathcal{S}_i$, the
algebras $\mathcal{J}_i$ and the homomorphisms $\isoName_i$ given $\mathcal{S}$~\cite[Chapter
6]{permenterThesis2017}. %These algorithms build on results of~\citet{jacobson1954structure}.

\subsubsection{Connections with *-subalgebras}\label{sec:starDecomp}
In some cases,  *-algebra techniques (currently used in the SDP literature) can find the decomposition of a Jordan
subalgebra $\mathcal{S}$ into its simple ideals, a crucial step in finding~\eqref{eq:decompForm}.  To explain, view $\mathbb{R}^{n \times n}$ as a *-algebra
with matrix multiplication as a product and transposition as a *-involution, and
let $\mathcal{M}\subseteq \mathbb{R}^{n \times n}$ denote the *-subalgebra generated by $\mathcal{S}$.
The \emph{Wedderburn decomposition}~\cite{wedderburn1908hypercomplex} of $\mathcal{M}$ is its direct-sum decomposition 
$\mathcal{M} =  \oplus^q_{i=1} \mathcal{M}_{i}$ into simple ideals $\mathcal{M}_i$. If $\mathcal{S} = \mathcal{M} \cap \mathbb{S}^n$, then the ideals of $\mathcal{M}$ identify the ideals of $\mathcal{S}$. Formally:
\begin{prop}[Ideals from the Wedderburn decomposition]
  Let $\mathcal{M}$ be the *-subalgebra of $\mathbb{R}^{n \times n}$ generated
  by a Jordan subalgebra $\mathcal{S}$ of $\mathbb{S}^n$. Let $\mathcal{M}$
  have Wedderburn decomposition $\mathcal{M} = \oplus^q_{i=1} \mathcal{M}_i$. If
  $\mathcal{S} = \mathcal{M} \cap \mathbb{S}^n$, then  $\oplus^q_{i=1} (\mathcal{M}_i \cap
  \mathbb{S}^n)$ is the decomposition of $\mathcal{S}$ into simple ideals.
  \begin{proof}
    We  need to show that $\mathcal{S} = \oplus^q_{i=1} (\mathcal{M}_i \cap \mathbb{S}^n)$
    and that $\mathcal{M}_i \cap \mathbb{S}^n$ is a simple ideal.
    
    To begin, write $X \in \mathcal{S}$ as $X=\sum^q_{i=1} X_i$ for $X_i \in \mathcal{M}_i$.  Then,
    $X=\sum^q_{i=1}\frac{1}{2}(X_i + X^T_i)$, where $X_i+X^T_i \in \mathcal{M}_i \cap
    \mathbb{S}^n$ since $\mathcal{M}_i$ is closed under transposition. Hence, $\mathcal{S} \subseteq \oplus^q_{i=1} (\mathcal{M}_i \cap \mathbb{S}^n)$. The reverse containment
    follows because $\mathcal{S} = \mathcal{M} \cap \mathbb{S}^n = (\oplus^q_{i=1} \mathcal{M}_i) \cap \mathbb{S}^n \supseteq \oplus^q_{i=1} (\mathcal{M}_i \cap \mathbb{S}^n)$.
   
    That $\mathcal{S}_i:=\mathcal{M}_i \cap \mathbb{S}^n$ is an ideal of $\mathcal{S}$ is obvious: if $X \in \mathcal{S}$ and $Y \in \mathcal{S}_i$ then $XY+YX \in \mathcal{M}_i$ since $\mathcal{M}_i$ is an ideal of $\mathcal{M}$,
    hence $\frac{1}{2}(XY+YX) \in \mathcal{S}_i$.
    Further, by the Artin-Wedderburn theorem,  each $\mathcal{M}_i$ is isomorphic to the *-algebra of real, complex, or quaternion
    matrices of some order; hence, $\mathcal{S}_i$ is isomorphic to the Hermitian matrices of real, complex, or quaternion entries of some order  and is therefore simple.
  \end{proof}
\end{prop}
Algorithms for finding the Wedderburn decomposition of *-subalgebras of
$\mathbb{R}^{n \times n}$
include~\cite{maehara2010numerical,eberly1996efficient}; see
also~\cite{gijswijt2010matrix,de2011numerical} for decompositions of complex *-algebras.

\begin{rem}
  A subalgebra $\mathcal{S}$ that satisfies $\mathcal{S} = \mathcal{M} \cap \mathbb{S}^n$ is called \emph{reversible}.
If a subalgebra is not reversible, one of its simple ideals is isomorphic to a spin-factor algebra.
Conversely, if a subalgebra is isomorphic to a spin-factor of dimension larger than 5, it is not reversible~\cite[Theorem 6.2.5]{hanche1984jordan}.
\end{rem}

\section{Minimum rank projections and admissible subspaces}\label{sec:findingsub}
  We now show how to find a projection $P_{\mathcal{S}} : \mathbb{S}^n \rightarrow \mathbb{S}^n$  satisfying  the
  Constraint Set Invariance  and  Unitality Conditions
  (Definitions~\ref{cond:CI}~and~\ref{defn:unital}), which, as argued in the
  previous section, allows one to reformulate the primal-dual
  pair~\eqref{sdp:main} over a symmetric cone isomorphic to $\mathbb{S}^n_{+}
  \cap \mathcal{S}$.  As we'll show, among projections that satisfy these
  conditions, there exists a unique one of minimum rank. Further, 
  a simple algorithm finds this projection for any instance of the primal-dual pair~\eqref{sdp:main}. This will follow by characterizing subspaces whose orthogonal
  projections satisfy these conditions. We define such a subspace as
  \emph{admissible}:
\begin{defn}\label{defn:admiss}
A subspace $\mathcal{S}$ is \emph{admissible} if its orthogonal projection $P_{\mathcal{S}} : \mathbb{S}^n \rightarrow \mathbb{S}^n$ satisfies  the Constraint Set Invariance and Unitality Conditions  (Definitions~\ref{cond:CI}-\ref{defn:unital}).
\end{defn}

Theorem~\ref{thm:posjor} provided a partial characterization of admissibility, showing that the ranges of positive, unital projections are the Jordan subalgebras of $\mathbb{S}^n$. To complete a characterization, we need the following result on invariance of the primal-dual affine sets.

\begin{lem}\label{lem:subspacechar}
	For affine sets $Y + \mathcal{L}$  and $C + \mathcal{L}^{\perp}$, let 
	$Y_{\mathcal{L}^{\perp}} \in \algName$ and $C_{\mathcal{L}}\in \algName$ denote the projections of $Y  \in \algName$ and $C \in \algName$ onto the subspaces    $\mathcal{L}^{\perp}$ and $\mathcal{L}$, respectively. The following are equivalent.

	\begin{enumerate}
		\item $C + \mathcal{L}$  and $Y + \mathcal{L}^{\perp}$ are invariant under the orthogonal projection $P_{\mathcal{S}} : \algName \rightarrow \algName$.
		\item\label{affine:statement:two} The subspace $\mathcal{S}$ contains $C_{\mathcal{L}}$ and $Y_{\mathcal{L}^{\perp}}$ and is  an invariant subspace of $P_{\mathcal{L}} : \algName \rightarrow \algName$, 
		i.e., $\mathcal{S}$ contains $Y_{\mathcal{L}^{\perp}}, C_{\mathcal{L}}$  and  $P_{\mathcal{L}}( \mathcal{S})$.
		
	\end{enumerate}
	\begin{proof}
		See  appendix (Section~\ref{sec:invarAffine}).
	\end{proof}
\end{lem}
 \begin{rem}The   invariant subspaces of $P_{\mathcal{L}}$  are precisely the
 invariant subspaces of $P_{\mathcal{L}^{\perp}}$~\cite[Proposition
 3.8]{farenick2012algebras}, which explains the asymmetry of
 statement Lemma~\ref{lem:subspacechar}-\eqref{affine:statement:two} with respect to
 $\mathcal{L}$ and $\mathcal{L}^{\perp}$. 
 \end{rem}
\noindent We also use the following well-known characterization of subalgebras, which follows given that $XY+YX = {(X+Y)}^2-X^2-Y^2$.
 \begin{lem}\label{lem:algchar}
   A subspace $\mathcal{S} \subseteq \mathbb{S}^n$ is a Jordan subalgebra of $\mathbb{S}^n$ (with product $\frac{1}{2} (XY+YX)$) if and only if $\mathcal{S} \supseteq \{ X^2 : X \in \mathcal{S} \}$.
 \end{lem} 
\noindent Combining Theorem~\ref{thm:posjor} with Lemmas~\ref{lem:subspacechar}-\ref{lem:algchar} yields our characterization of admissibility.
\begin{thm}\label{thm:admisschar}
	A subspace $\mathcal{S}$ is admissible (Definition~\ref{defn:admiss}) if and only if it satisfies the following conditions: 
	\begin{align*} 
	\mathcal{S} &\ni Y_{\mathcal{L}^{\perp}}, C_{\mathcal{L}} ,    \\ 
	\mathcal{S} &\supseteq P_{\mathcal{L}}(\mathcal{S}),   \\
	\mathcal{S} &\supseteq \{ X^2 : X \in \mathcal{S} \},
	\end{align*}
	where $Y_{\mathcal{L}^{\perp}}$  and $C_{\mathcal{L}}$ are as in Lemma~\ref{lem:subspacechar} and $\mathcal{L}$ is the linear subspace of the primal-dual pair~\eqref{sdp:main}.
\end{thm}

\subsection{Optimal subspaces and minimum rank projections}
\noindent Theorem~\ref{thm:admisschar} shows that the (arbitrary) intersection of admissible subspaces is admissible.  This motivates the following definition.
\begin{defn}\label{defn:optsub}
	The \emph{optimal admissible subspace} $\mathcal{S}_{opt}$ is the intersection of all admissible subspaces:
	\[
	\mathcal{S}_{opt} = \bigcap \{ \mathcal{S} : \mbox{ $\mathcal{S}$  is admissible} \}.
	\]
\end{defn}
\noindent Admissibility of $\mathcal{S}_{opt}$ yields the following corollary of Theorem~\ref{thm:admisschar}.
\begin{cor}
  The map $P_{\mathcal{S}_{opt}} : \mathbb{S}^n \rightarrow \mathbb{S}^n$ is the minimum-rank orthogonal projection satisfying the Constraint Set Invariance and Unitality Conditions (Definitions~\ref{cond:CI}-\ref{defn:unital}).
\end{cor}

\subsection{Solution algorithm}
 Theorem~\ref{thm:admisschar} also suggests a procedure for finding $\mathcal{S}_{opt}$.
 First, initialize $\mathcal{S}$ to the subspace spanned by $C_{\mathcal{L}}$ and $Y_{\mathcal{L}^{\perp}}$. Then, add $P_{\mathcal{L}}(\mathcal{S})$ and the span of $\{ X^2 : X \in \mathcal{S} \}$
to $\mathcal{S}$ in an alternating fashion, terminating when the resulting ascending chain of subspaces  stabilizes. Formally:
\begin{thm}\label{thm:subalg}
	The  optimal admissible subspace $\mathcal{S}_{opt}$ (Definition~\ref{defn:optsub}) is the output of the following algorithm:
	\begin{center} \qquad \qquad
		\DontPrintSemicolon
		\begin{minipage}{.35 \textwidth}
			$\mathcal{S} \leftarrow \Span \{ C_{\mathcal{L}}, Y_{\mathcal{L}^{\perp}  }\} $ \\
                        \Repeat{\rm ascending chain stabilizes}{
				$\mathcal{S} \leftarrow \mathcal{S} + P_{\mathcal{L}}(\mathcal{S})$ \\
                                $\mathcal{S} \leftarrow \mathcal{S}+ \Span \{ X^2 : X \in \mathcal{S} \}$
			}  
		\end{minipage}
	\end{center} 
	where   $C_{\mathcal{L}}, Y_{\mathcal{L}^{\perp}}$ are as in Lemma~\ref{lem:subspacechar} and $\mathcal{L}$ is the linear subspace of the primal-dual pair~\eqref{sdp:main}.
	\begin{proof}
		The algorithm computes an ascending chain of subspaces in finite dimensions 
		which must stabilize to a subspace $\mathcal{\hat S}$. Stabilization implies that
		\[
		\mathcal{\hat S} = \mathcal{\hat S} + P_{\mathcal{L}}(\mathcal{\hat S}), \qquad \mathcal{\hat S} = \mathcal{\hat S}+ \Span \{ X^2 : X \in \mathcal{\hat S} \}
		\]
		which, since $\mathcal{\hat S} \ni  C_{\mathcal{L}}, Y_{\mathcal{L}^{\perp}}$, shows that $\mathcal{\hat S}$ is admissible (Theorem~\ref{thm:admisschar}).  At every iteration, $\mathcal{S}$ is a subset of  $\mathcal{S}_{opt}$ (by induction); hence, $\mathcal{\hat S}$ is a subset. But since $\mathcal{\hat S}$ is admissible, it contains $\mathcal{S}_{opt}$ (Definition~\ref{defn:optsub}).  We conclude that $\mathcal{\hat S} = \mathcal{S}_{opt}$.
	\end{proof}
\end{thm}
\noindent 
Note executing this algorithm may be impractical if storing a basis for
$\mathcal{S}$ is impractical.  To deal with this, we introduce variations in
Section~\ref{sec:comb} that restrict to subspaces spanned by bases with
efficient combinatorial representations.

\section{Optimal decompositions}\label{sec:optdirectsum}

Admissible subspaces are necessarily  subalgebras of $\algName$
(Theorem~\ref{thm:admisschar}). As a consequence,  each admissible subspace has an orthogonal direct-sum decomposition
into simple ideals (Section~\ref{sec:subalgstruc}). We now prove 
the direct-sum decomposition of $\mathcal{S}_{opt}$ is optimal in a precise sense.

Our notion of optimality is in terms of the \emph{rank vector} of an algebra $\mathcal{W}=\oplus^w_{i=1} \mathcal{W}_i$
\[
r_{\mathcal{W}} := (\rank \mathcal{W}_1,\rank \mathcal{W}_2,\ldots,\rank \mathcal{W}_w),
\] 
where each $\mathcal{W}_i$ is a simple ideal and   $\rank \mathcal{W}_i$ is the maximum  number of distinct eigenvalues
of an $X \in \mathcal{W}_i$; as an example, the rank vectors of $\mathbb{S}^{n_1} \oplus \mathbb{S}^{n_2}$  and $\mathbb{R}^4_{+}$ are $(n_1,n_2)$ and $(1,1,1,1)$, respectively.
Specifically, we show that the rank vector of   $\mathcal{S}_{opt}$ is \emph{weakly majorized} by that of any other
admissible  subspace. Among other things, this means that $\mathcal{S}_{opt}$ minimizes the rank vector's largest element and the sum of its elements.
\begin{defn}\label{defn:majorize}
	The vector $x  \in \mathbb{Z}^m$   \emph{weakly majorizes}  $y  \in \mathbb{Z}^n$ if
	\[
          \sum^{\min{\{\ell,m\}}}   _{i=1} {[x^{\downarrow} ]}_i \ge \sum^{\min{\{\ell,n\}}}_{i=1} {[y^{\downarrow} ]}_i \qquad \forall \ell \in \{1,\ldots,\max{\{m,n\}}     \},
	\]
	where $x^{\downarrow} $ and $y^{\downarrow} $ denote  $x $ and $y $ with entries sorted in descending order.
\end{defn}

\begin{ex}
For the following   subalgebras $\mathcal{U}_i$ (each parametrized by  $t \in \mathbb{R}^m$), the rank vector  $r_{\mathcal{U}_i}$ weakly majorizes $r_{\mathcal{U}_{i+1}}$:
\begin{align*} 
&\begin{array}{cc}
\mathcal{U}_1 :=& \left(\begin{array}{ccccc} t_1 & t_2 & 0 & 0 & 0 \\ t_2 & t_3 & 0 & 0 & 0 \\ 0 & 0 & t_4 & t_5 & t_6 \\ 0 & 0 & t_5 & t_7 & t_8 
\\ 0 & 0 & t_6 & t_8 & t_9 
\end{array}\right)  \\\\
&r_{\mathcal{U}_1} = (2,3) 
\end{array}
\qquad
\begin{array}{cc}
\mathcal{U}_2 :=& \left(\begin{array}{ccccc} t_1 & t_2 & 0 & 0 & 0 \\ t_2 & t_3 & 0 & 0 & 0 \\ 0 & 0 & t_4 & t_5 & 0 \\ 0 & 0 & t_5 & t_7 & 0 
\\ 0 & 0 & 0 & 0 & t_9 
\end{array}\right)  \\\\
&r_{\mathcal{U}_2} = (2,2,1) 
\end{array} \\ \\
&\begin{array}{cc}
\mathcal{U}_3 :=& \left(\begin{array}{ccccc} t_1 & t_2 & 0 & 0 & 0 \\ t_2 & t_3 & 0 & 0 & 0 \\ 0 & 0 & t_1 & t_2 & 0 \\ 0 & 0 & t_2 & t_3 & 0 
\\ 0 & 0 & 0 & 0 & t_4 
\end{array}\right)  \\\\
&r_{\mathcal{U}_3} = (2,1) 
\end{array}
\qquad 
\begin{array}{cc}
\mathcal{U}_4 :=& \left(\begin{array}{ccccc} t_1 & 0 & 0 & 0 & 0 \\ 0 & t_1 & 0 & 0 & 0 \\ 0 & 0 & t_1 & 0 & 0 \\ 0 & 0 &  0 & t_1 & 0 
\\ 0 & 0 & 0 & 0 & t_2 
\end{array}\right)  \\\\
&r_{\mathcal{U}_4   } = (1,1) 
\end{array} 
\end{align*} 
%Also of note are the subalgebras $\mathcal{U}_3$ and $\mathcal{U}_4$; despite having \emph{three} nonzero blocks,  $\mathcal{U}_3$ is isomorphic to a product  of \emph{two} simple algebras since its second two-by-two block is a copy of the first; similar remarks apply to $\mathcal{U}_4$.
\end{ex}

\noindent We now state our result.
\begin{thm}\label{thm:OptimalDecomp}
         Let $\mathcal{W} \subseteq \mathbb{S}^n$  be any admissible   subspace
         (Definition~\ref{defn:admiss}).   Let the optimal admissible  subspace
         $\mathcal{S}_{opt}$ (Definition~\ref{defn:optsub}) and $\mathcal{W}$
         have the following decompositions into simple ideals:
	\[
	\mathcal{S}_{opt}=\oplus^s_{i=1}\mathcal{S}_i, \qquad \mathcal{W}= \oplus^w_{k=1} \mathcal{W}_k.
	\]
        Then, $r_\mathcal{W} := (\rank \mathcal{W}_1,\ldots,\rank \mathcal{W}_w)$ weakly majorizes $r_{\mathcal{S}_{opt}}:= (\rank \mathcal{S}_1,\ldots,\rank \mathcal{S}_s)$.

\end{thm}  

\noindent To prove this theorem, we will only use the fact that $\mathcal{S}_{opt}$ is a subalgebra of all other  admissible subspaces, which is immediate from its definition and Theorem~\ref{thm:admisschar}.

\subsection{Proof of Theorem~\ref{thm:OptimalDecomp}}
\renewcommand{\algName}{\mathcal{J}}
We prove the theorem by showing a  more general result (Theorem~\ref{thm:rankSubAlg}) about the rank vectors of Euclidean Jordan
algebras and their subalgebras. To our knowledge, these results are new.
Towards this, we let $x \circ y$ denote the Jordan product of an abstract
Euclidean Jordan algebra $\mathcal{J}$ and recall some standard definitions. An
\emph{idempotent} is an $x \in \algName$ satisfying $x \circ x = x$. An
idempotent is \emph{primitive} if it is nonzero and doesn't equal the sum  of
two different nonzero idempotents. Finally, as mentioned in
Section~\ref{sec:subalgstruc}, an algebra $\mathcal{J}$ is
simple if its only ideals are $\mathcal{J}$ and $\{0\}$. We start with a needed
technical lemma.
\begin{lem}\label{lem:jorprojprop}

	Let $\algName$  be a Euclidean Jordan algebra and let $\mathcal{V} \subseteq \algName$ be a subalgebra that  is simple (viewed as an algebra).  Let $\algName = \oplus^w_{k=1} \algName_k$ denote the orthogonal direct-sum decomposition of $\algName$ into simple ideals.  Finally, let $\Phi_k : \algName \rightarrow \algName$ denote the orthogonal projection onto $\algName_k$.  The following statements hold for all $k \in [w]$, where $[w] := \{1,\ldots,w\}$:
	
	\begin{enumerate}
		\item If $e \in \algName$ is an idempotent, then $\Phi_k e $ is an idempotent.
                \item    If  $e,f \in \algName$ are idempotents and $\langle e , f \rangle = 0$, then $\langle \Phi_{k} e, \Phi_{k}f \rangle = 0$.
		\item Suppose $e,f \in \mathcal{V}$ are nonzero idempotents.  If $\Phi_k e \ne  0$, then $\Phi_k f \ne 0$.
	\end{enumerate}
	
	\begin{proof}
          Since $\algName_k$ is a simple ideal, the projection map $\Phi_k$ from $\algName$ onto $\algName_k$  is a Jordan homomorphism by~\cite[Lemma 2.5.6]{hanche1984jordan}; hence, $\Phi_k e \circ \Phi_k e = \Phi_k (e^2) =  \Phi_k e$,
		showing the first statement.
		
		For the second statement, recall $\algName = \oplus^w_{k=1} \algName_k$ is an orthogonal direct-sum decomposition of $\algName$. We conclude 
		\[
		e = \sum^w_{k=1} \Phi_k e, \qquad  f = \sum^w_{k=1} \Phi_k e.
		\]
                Since $\langle \Phi_i e  , \Phi_j f  \rangle \ge 0$ (since the cone-of-squares is self-dual) and
		\[
		\langle e , f \rangle = \sum^w_{i=1} \sum^w_{j=1} \langle \Phi_i e  , \Phi_j f  \rangle,
		\]
		$\langle \Phi_i e  , \Phi_j f  \rangle = 0$ if $\langle e , f \rangle =0$.

		For the third statement, view $\mathcal{V}$ as a simple algebra and let  $e=\sum^q_{i=1} e_i$ and $f=\sum^r_{j=1} f_j$ denote the decompositions of $e$ and $f$ into primitive idempotents of $\mathcal{V}$. Then, there exists $t \in \mathcal{V}$ (depending on $i$ and $j$) such that $e_i= 2 t \circ (t \circ f_j)- t^2 \circ f_j$~\cite[Corollary IV.2.4]{faraut1994analysis}.  Since $\Phi_k$
		is a homomorphism, 
		\begin{align*}
		\Phi_k e_i   &= \Phi_k( 2 t \circ (t \circ f_j)- t^2 \circ f_j     ) \\&=   \Phi_k(2t) \circ ( \Phi_k t  \circ \Phi_k f_j )- \Phi_k t^2 \circ \Phi_k f_j        
		\end{align*}
		showing $\Phi_k f_j   \ne 0$ if  $\Phi_k e_i  \ne 0$.  Since
		\[
		\Phi_k e  = \sum^q_{i=1}  \Phi_k e_i  , \;\;\; \Phi_k f  = \sum^r_{j=1}  \Phi_k f_j ,
		\]
		and $\Phi_k e_i $ and $\Phi_k f_j $ are idempotents and hence in the cone-of-squares, it follows  $\Phi_k f   \ne 0$ if $\Phi_k e  \ne 0$.
		
	\end{proof}
\end{lem}
\noindent The mentioned results on rank vectors and subalgebras follow.
\begin{thm}[Subalgebras and rank vectors]\label{thm:rankSubAlg}
	Let $\mathcal{S}=\oplus^s_{i=1} \mathcal{S}_i$ and $\mathcal{W} = \oplus^w_{k=1} \mathcal{W}_k$  be  Jordan subalgebras of $\algName$, where $\mathcal{S}_i$ and $\mathcal{W}_k$ are simple ideals of $\mathcal{S}$ and $\mathcal{W}$ (viewed as algebras), respectively. Suppose $\mathcal{S} \subseteq \mathcal{W}$. The following statements hold:
	
	\begin{enumerate}
          \item  For each $k \in [w]$, let  $I_k := \{ i \in [s] : \mathcal{S}_i \not\subseteq {(\mathcal{W}_k)}^{\perp} \}$. Then, for all $k \in [w]$,
		\[
		\rank \mathcal{W}_k \ge \sum_{i \in I_k} \rank \mathcal{S}_i.
		\]
		
		\item  The vector $r_\mathcal{W}$ weakly majorizes $r_\mathcal{S}$, where
		\[
		r_\mathcal{W}:= (\rank \mathcal{W}_1,\ldots,\rank \mathcal{W}_w), \qquad r_\mathcal{S}:= (\rank \mathcal{S}_1,\ldots,\rank \mathcal{S}_s).
		\]
		
	\end{enumerate}

	\begin{proof}  
          First note $\mathcal{S}_i$  contains a set $\mathcal{E}_i := {\{e^i_j\}}^{  \rank \mathcal{S}_i}_{j=1}$ of  pairwise-orthogonal  idempotents. Further, if $i \in I_k$, then $\Phi_k e \ne 0$
		for a nonzero idempotent  $e$ in  $\mathcal{S}_i$.      We conclude all elements of $\{ \Phi_k f : f \in \cup_{i \in I_k} \mathcal{E}_i\}$  are nonzero  (Lemma~\ref{lem:jorprojprop}-3); moreover, they are idempotent (Lemma~\ref{lem:jorprojprop}-1) and pairwise orthogonal (Lemma~\ref{lem:jorprojprop}-2).  It follows $\mathcal{W}_k$ contains
		at least $\sum_{i \in I_k} \rank \mathcal{S}_i$ nonzero idempotents that are pairwise orthogonal. Hence,  $\rank \mathcal{W}_k \ge \sum_{i \in I_k} \rank \mathcal{S}_i$.

		For the second statement, we note the first implies the following: for each $\ell \in \max\{s,w\}$, there is a subset $T \subseteq [w]$  for which 
		\[
                  \sum_{k \in T} \rank \mathcal{W}_k \ge \sum_{k \in T} \sum_{i \in I_k} \rank \mathcal{S}_i \ge \sum^{\min{\{\ell,s\}}}_{i=1} {[r^{\downarrow}_\mathcal{S}]}_i.
		\]
                Specifically, letting $\pi$ be a permutation of $[s]$ satisfying ${[r^{\downarrow}_\mathcal{S}]}_i = {[r_\mathcal{S}]}_{\pi(i)}$,  we can choose $T$ to be subset of $[w]$ that satisfies
		\[
                  \cup_{k \in T} I_k \supseteq  {\{ \pi(i) \}}^{\min \{\ell,s\}}_{i=1}.
		\] 
		Further, we can choose  $T$ to have  $|T| \le \min{\{\ell,w\}}$, which implies that 
		\[
                  \sum^{\min{\{\ell,w\}}}_{i=1} {[r^{\downarrow}_\mathcal{W}]}_i \ge \sum_{k \in T} \rank \mathcal{W}_k.
		\]
                Hence,  the majorization inequality $\sum^{\min{\{\ell,w\}}}_{i=1} {[r^{\downarrow}_\mathcal{W}]}_i \ge \sum^{\min{\{\ell,s\}}}_{i=1} {[r^{\downarrow}_\mathcal{S}]}_i$ holds.
		
		%\fbox{Definition of $T$}:
		
		%	\[
		%	T = \cup_i \{ k : \pi(i) \in I_k  \}
		%	\]
		
	%	For the second statement,   Consider the bipartite graph $\mathcal{G}$ with node sets
	%	$[s]$ and $[w]$, where $(i,k)$ are adjacent iff $i \in I_k$. 	Let $\pi$ be a permutation of $[s]$ satisfying $[r^{\downarrow}_\mathcal{S}]_i = [r_\mathcal{S}]_{\pi(i)}$.

	%	 Now, fix $\ell \in \max\{s,w\}$ and let $T \subseteq [w]$ be
		%a minimum cardinality subset of $[w]$ such that each node of $\{ \sigma(i)  \}^{\min\{\ell,s\} }_{i=1}$ is adjacent to at least one node of $T$ in the graph  $\mathcal{G}$. Then, 
		%	\[
		%	\sum_{k \in T} \rank \mathcal{W}_k \ge \sum_{k \in T} \sum_{i \in I_k} \rank \mathcal{S}_i \ge \sum^{\min{\{l,s\}}}_{i=1} [r^{\downarrow}_\mathcal{S}]_i.
		%	\]

	\end{proof}
\end{thm}
We see that Theorem~\ref{thm:OptimalDecomp} follows immediately from the second statement of Theorem~\ref{thm:rankSubAlg} since, as mentioned,
$\mathcal{S}_{opt}$ is a subalgebra of all other  admissible subspaces (Definition~\ref{defn:majorize}). 
\renewcommand{\algName}{\mathbb{S^n}}

\section{Combinatorial variations}\label{sec:comb}

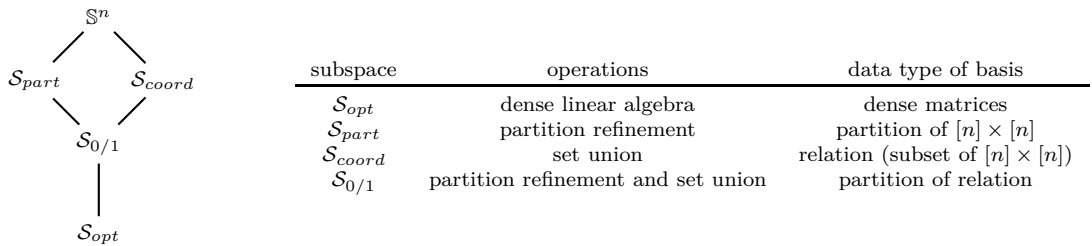
\begin{figure}
	
	\begin{minipage}{0.3\textwidth}
		\centering
		\begin{tikzpicture}
		[align=center,node distance=1.2cm]
		\node (full) at (0,0) {$\mathbb{S}^n$};
		\node [below left of=full] (part)  {$\mathcal{S}_{part}$};
		\node [below right of=full] (coord)  {$\mathcal{S}_{coord}$};  
		\node [below left of=coord] (zeroone)  {$\mathcal{S}_{0/1}$};
		\node [below of=zeroone] (opt)  {$\mathcal{S}_{opt}$};
		\draw [thick] (full) -- (coord);
		\draw [thick] (full) -- (part);
		\draw [thick] (coord) -- (zeroone);
		\draw [thick] (part) -- (zeroone);
		\draw [thick] (zeroone) -- (opt);
		\end{tikzpicture}
	\end{minipage} 
	\begin{minipage}{0.45\textwidth}
		\centering
		\begin{tabular}{ccc}
			subspace &  operations & data type of basis\\ \toprule
			$\mathcal{S}_{opt}$ & dense linear algebra & dense matrices \\  
			$\mathcal{S}_{part}$ & partition refinement & partition of $[n] \times [n]$ \\ 
			$\mathcal{S}_{coord}$ & set union & relation (subset of $[n] \times [n]$) \\ 
			$\mathcal{S}_{0/1}$ & partition refinement and set union  &  partition of  relation \\ 
		\end{tabular}
	\end{minipage}
        \caption{Hasse diagram of subspace inclusions,     key algorithmic operation needed to find subspace,  and data type (mathematical object) used to represent a basis.}\label{fig:hasse}
\end{figure}

This section introduces  combinatorial restrictions on admissible subspaces
(Definition~\ref{defn:admiss}), aiming to reduce the cost of storing a basis.
We consider three types of subspaces (Figure~\ref{fig:hasse}). The first two
types have bases encoded by \emph{relations} and \emph{partitions}, respectively. The third type is a
common generalization, whose discussion we defer to the end of this section. 

To begin, let $[n] = \{1,\ldots,n\}$. For a relation  $\mathcal{R} \subseteq [n] \times [n]$,
let $\mathcal{B}_{\mathcal{R}}:= \left\{ E_{ij} + E_{ji}   : (i,j) \in \mathcal{R}   \right\}$, where $E_{ij}$ is
the standard basis matrix of $\mathbb{R}^{n \times n}$ nonzero (and equal to 1) only at its ${(i,j)}$-th entry. For a partition $\mathcal{P}$ of $[n] \times [n]$, let $\mathcal{B}_{\mathcal{P}}$ denote the corresponding set of characteristic matrices---i.e., let $\mathcal{B}_{\mathcal{P}} :=  \{  \sum_{(i,j) \in C} E_{ ij }  :  C \in \mathcal{P}\}$ where  $C \subseteq [n] \times [n]$ denotes a subset in $\mathcal{P}$.

\begin{defn}
	A \emph{coordinate subspace}  is the span of $\mathcal{B}_{\mathcal{R}}$ for some relation $\mathcal{R} \subseteq [n] \times [n]$.
	A \emph{partition subspace}  is the span of $\mathcal{B}_{\mathcal{P}}$ for some partition $\mathcal{P} \subseteq [n] \times [n]$.
\end{defn}

 \begin{ex}
 	The following subspaces $\mathcal{S}_1$ and $\mathcal{S}_2$ are coordinate and partition subspaces, respectively.
 	\[
 	\mathcal{S}_1  = \left\{  \begin{bmatrix}
 	a & b & 0 \\ 
 	b & c & d \\ 
 	0 & d & e
 	\end{bmatrix} : (a,b,c,d) \in \mathbb{R}^4  \right\} \qquad
 	\mathcal{S}_2 = \left\{  \begin{bmatrix}
 	a & a & b \\ 
 	a & a & b \\ 
 	b & b & c 
 	\end{bmatrix} : (a,b,c) \in \mathbb{R}^3 \right\}.
 	\]
        Specifically, $\mathcal{S}_1$ equals the span of $\mathcal{B}_{\mathcal{R}}$ for the relation
	\[
	  \mathcal{R} = \left\{ (1,1), (1,2),(2,1), (2,2),(2,3),(3,2),(3,3)\right\},
	 \]
	 and $\mathcal{S}_2$ equals the span  of $\mathcal{B}_{\mathcal{P}}$ for the partition 
	\[
	 \mathcal{P} = \bigg\{  \{ (1,1),(1,2),(2,1),(2,2)     \}, \;\;  \{ (1,3),(2,3), (3,1), (3,2) \}, \;\; \{ (3,3)     \} \bigg\}.
	\]
 \end{ex}
\noindent We seek the following variants of $\mathcal{S}_{opt}$:
\begin{align*}
  \mathcal{S}_{coord} &:= \bigcap \{  \mathcal{S} \subseteq \mathbb{S}^n : \mbox{$\mathcal{S}$ is admissible and a coordinate subspace}\},\\
  \mathcal{S}_{part} &:=  \bigcap \{  \mathcal{S} \subseteq \mathbb{S}^n : \mbox{$\mathcal{S}$ is admissible and a partition subspace}\}.
\end{align*}
The families of admissible, coordinate, and partition subspaces are all closed under intersection.
Hence, $\mathcal{S}_{coord}$ is both admissible and a coordinate subspace.   Similar remarks apply for $\mathcal{S}_{part}$.

Though coordinate subspaces are a highly restricted family, our conference paper~\cite{permenterCoord2015}
illustrates $\mathcal{S}_{coord}$ can have small dimension for SDPs arising in polynomial optimization.
Partition subspaces also arise naturally in symmetry reduction.  Indeed, the fixed-point subspace (Section~\ref{sec:suffgen}) 
\[
  \mathcal{M}_\mathcal{G} = \{ X \in \mathbb{R}^{n \times n} : P X P^T =  X \;\; \forall P \in \mathcal{G} \}
\]
is a partition subspace (of $\mathbb{R}^{n \times n}$) and
$\mathcal{M}_{\mathcal{G}} \cap \mathbb{S}^n$ a partition subspace of
$\mathbb{S}^n$ when $\mathcal{G}$ is a group of permutation matrices.  The
partition $\mathcal{P}$ of $[n] \times [n]$ that induces $\mathcal{M}_{\mathcal{G}}$ arises from
the \emph{orbits}  $\{P E_{ij} P^T : P \in \mathcal{G}\}$ of the standard basis
matrices $E_{ij}\in \mathbb{R}^{n \times n}$; precisely,  $(i,j)$ and $(k, l)$
are in the same class of $\mathcal{P}$ if the orbit $\{ P E_{ij} P^T : P \in \mathcal{G} \}$
contains  $E_{kl}$.   (Such a partition  is called a  \emph{Schurian coherent
configuration}~\cite{higman1987coherent}.)

%A partition $\mathcal{P}$ is called
%a coherent configuration if the span of $\mathcal{B}_{\mathcal{P}}$ is a *-subalgebra
%of $\mathbb{R}^{n \times n}$ that contains $I$, which contains all subspaces of the form $\mathcal{C}_{\mathcal{G}}$.

%Based on this, let us define \emph{coherent Jordan algebras} as the Jordan subalgebras of $\mathbb{S}^n$ that are also partition subspaces.  To our knowledge, it is unknown if a converse of this proposition holds, i.e., it is unknown if all coherent Jordan algebras can be written as the symmetric part of a coherent algebra.~\citeauthor{cameron2003coherent} has posed a similar open question~\cite[page 8]{cameron2003coherent}.

% specifically, it indicates if they always have a sparse basis and a sparse \emph{decomposition}.  By sparse decomposition, we mean their minimal ideals (when viewed as Jordan algebras)  have a sparse basis.  Note a sparse basis is still valuable even if the subspace does not have a sparse decomposition, as we illustrate in the example section.

%A $0/1$ basis of the first type is the set of  $0/1$ characteristic matrices of  a partition.  A  basis of the second type spans a coordinate subspace  of $\mathbb{S}^n$---that is, it spans a subspace specified only by a sparsity pattern.  A  basis of the third type is a simply a set of $0/1$ matrices that are pairwise orthogonal.  Note this first type closely relates to coherent-configurations~\cite{cameron2003coherent},   the second type   the topic of our conference paper~\cite{permenterCoord2015}, and the third type a generalization of the first two. 

%Specifically, three alternative methods are given for optimizing over  subspaces  with  different   combinatorial descriptions.  We first optimize over subspaces spanned by the characteristic matrices of  partitions, which is inspired by coherent-configuration-based reduction methods~\cite{de2010exploiting}. We then optimize over the  coordinate subspaces of $\mathbb{S}^n$; note we originally presented this method in the conference paper~\cite{permenterCoord2015}. Finally, we optimize over subspaces that have   an orthogonal basis of $0/1$ matrices---a set of subspaces that contains the first two.  In each case, we store a basis for $\mathcal{S}$ combinatorially; specifically, letting $n:= \{1,\ldots,n\}$,  we  first use  a partition  $\mathcal{P}$ of $[n] \times [n]$, next  a  relation $\mathcal{R} \subseteq [n]\times [n]$, and  third  a  partition of a relation.

% and involve partition refinement and set union as their basic operations, respectively.  An algorithm for the third, in a sense, combines these operations.

\subsection{Modified algorithms}\label{sec:algsectionComb}

To find $\mathcal{S}_{part}$ or $\mathcal{S}_{coord}$, we modify the
Theorem~\ref{thm:subalg} algorithm line-by-line to operate on partitions or
relations instead of subspaces. These modified algorithms 
first represent the image  of a coordinate/partition  subspace  $\mathcal{S}$ under the maps $X \mapsto X^2$
and $P_{\mathcal{L}} : \mathbb{S}^n \rightarrow \mathbb{S}^n$  with  a polynomial matrix, an idea inspired by~\cite{weisfeiler1977construction}; see
also \cite[Section 5]{babel2010algebraic}. They then refine/grow a
partition/relation based on the unique/nonzero entries of this polynomial matrix. These
algorithms are explicitly given in Figure~\ref{fig:algComb}.
They leverage the following notation (Definitions~\ref{defn:polymat}-\ref{defn:objects}).
\begin{defn}\label{defn:polymat}
	 For a finite set  $\mathcal{B} \subset \mathbb{S}^n$, let $f_{X^2}(\mathcal{B})$ and $f_{\mathcal{L}}(\mathcal{B})$ denote the  polynomial matrices
	 \[
           f_{X^2}(\mathcal{B}) := {\left(\sum_{B \in \mathcal{B}} t_{B} B \right)}^2 \qquad f_{\mathcal{L}}(\mathcal{B}) :=  \sum_{B \in \mathcal{B}} t_{B} P_{\mathcal{L}}(B) ,
	 \]
         \noindent where   $[t_{B}]_{B \in \mathcal{B}}$ is a vector of commuting\footnote{Note that the related algorithm~\cite[Section 5]{babel2010algebraic} uses \emph{noncommuting} indeterminates.} indeterminates indexed by $\mathcal{B}$. 
	
\end{defn}
\noindent If $\mathcal{S}$ is the span of $\mathcal{B}$, then the set of point evaluations of $f_{X^2}(\mathcal{B})$ equals $\{ X^2 : X \in \mathcal{S} \}$, i.e.,
 \[
   \{ X^2 : X \in \mathcal{S} \} = \{ {f_{X^2}(\mathcal{B})}_{|t_\mathcal{B}=t^*} : t^* \in \mathbb{R}^{|\mathcal{B} |} \},
 \]
and similarly for  $f_{\mathcal{L}}(\mathcal{B})$.  The following example illustrates this notation.
\begin{ex}
	For  $\mathcal{B} = \{U,V,W\}$, where

	\begin{align*}
	U = \left(\begin{array}{cccc} 0 & 0 & 1 & 0\\ 0 & 0 & 0 & 1\\ 1 & 0 & 0 & 0\\ 0 & 1 & 0 & 0 \end{array}\right), \;\; 
	V = \left(\begin{array}{cccc} 0 & 0 & 0 & 0\\ 0 & 0 & 0 & 0\\ 0 & 0 & 0 & 1\\ 0 & 0 & 1 & 0 \end{array}\right), \;\;
	W = \left(\begin{array}{cccc} 1 & 0 & 0 & 0\\ 0 & 1 & 0 & 0\\ 0 & 0 & 0 & 0\\ 0 & 0 & 0 & 0 \end{array}\right),
	\end{align*}
        we have $f_{X^2}(\mathcal{B}) := {(t_{U} U + t_{V} V + t_{W} W)}^2$.  Expanding (and using the identities $t_{U} t_{V} = t_{V} t_{U}$ and $t_{U} t_{W} = t_{W} t_{U}$)  shows that
	\[
	f_{X^2}(\mathcal{B}) =  \left(\begin{array}{cccc} 
	{t_U}^2 + {t_W}^2 & 0 & t_U t_W & t_U t_V \\ 
	0 & {t_U}^2 + {t_W}^2 & t_U t_V & t_U t_W\\ 
	t_U t_W & t_U t_V & {t_U}^2 + {t_V}^2 & 0\\ 
	t_U t_V & t_U t_W & 0 & {t_U}^2 + {t_V}^2 
	\end{array}\right).
	\]
\end{ex}
\begin{defn}\label{defn:objects} For an $n \times n$ polynomial matrix  $X$, let $\relation(X)$ denote the subset of $(i,j) \in [n] \times [n]$ for which $X_{ij}$ is not the zero polynomial.  Similarly, let $\partition(X)$ denote the partition of $[n] \times [n]$ induced
  by the unique polynomial entries of $X$, i.e., $(i,j)$ and $(k,l)$ are in the same class of $\partition(X)$ if and only if $X_{ij}$ and $X_{kl}$ are the same polynomial.
\end{defn}
\addtocounter{ex}{-1}
\begin{ex}[continued]
	For the polynomial matrix $f_{X^2}(\mathcal{B})$ of the previous example,
	the relation $\relation \left(f_{X^2}(\mathcal{B})\right)$ is the complement of $\{ (1,2),(2,1),(3,4),(4,3) \} \subseteq [n] \times [n]$ (where $n=4$.) The partition $\partition \left(f_{X^2}(\mathcal{B})\right)$  has  characteristic matrices  
	{\smaller
                \begin{align}\label{ex:Partition}
		\begin{array}{ccccc} 
		\left(\begin{array}{cccc} 
		1 & 0 & 0 & 0 \\ 
		0 & 1 & 0 & 0\\ 
		0 & 0 & 0 & 0 \\ 
		0 & 0 & 0 & 0
		\end{array}\right) & \left(\begin{array}{cccc} 
		0 & 1 & 0 & 0 \\ 
		1 & 0 & 0 & 0\\ 
		0 & 0 & 0 & 1 \\ 
		0 & 0 &  1 & 0
		\end{array}\right)  & \left(\begin{array}{cccc} 
		0 & 0 & 1  & 0 \\ 
		0 & 0 & 0 & 1\\ 
		1 & 0 & 0  & 0 \\ 
		0 & 1 & 0 & 0
		\end{array}\right)  & \left(\begin{array}{cccc} 
		0 & 0 & 0 & 1 \\ 
		0 & 0 & 1 & 0\\ 
		0 & 1 & 0 & 0 \\ 
		1 & 0 & 0 & 0 
		\end{array}\right) &  \left(\begin{array}{cccc} 
		0 & 0 & 0 & 0 \\ 
		0 & 0 & 0 & 0\\ 
		0 & 0 & 1 & 0 \\ 
		0 &0 & 0 & 1
		\end{array}\right), \\\\
		{t_U}^2 + {t_W}^2 & 0  & t_U t_W  & t_U t_V   & {t_U}^2 + {t_V}^2,
		\end{array}
		\end{align} 
	}
	where we've labeled each matrix by the associated polynomial entry of $f_{X^2}(\mathcal{B})$.
\end{ex}

\begin{figure} \centering
	\begin{subfigure}{.45\textwidth}
	\begin{algorithm}[H]
		$\mathcal{R} \leftarrow \relation \left( C_{\mathcal{L}} \right) \cup  \relation \left( Y_{\mathcal{L}^{\perp}} \right)$ 
		\\ \DontPrintSemicolon
                \Repeat{\rm ascending chain of relations $\mathcal{R}$ stabilizes.}{  
			$\mathcal{R} \leftarrow  \mathcal{R} \cup \relation \left( f_{\mathcal{L}}( \mathcal{B}_\mathcal{R}) \right)$  \\
			$\mathcal{R} \leftarrow  \mathcal{R} \cup \relation \left( f_{X^2} ( \mathcal{B}_\mathcal{R}) \right)$ 
		}   
	\end{algorithm}	
	\end{subfigure}% need this comment symbol to avoid overfull
	\begin{subfigure}{.45\textwidth}
		\begin{algorithm}[H]
			\DontPrintSemicolon
			$\mathcal{P}  \leftarrow  \partition (C_{\mathcal{L}})  \bigvee  \partition (Y_{\mathcal{L}^{\perp}})$ \\
                        \Repeat{\rm ascending chain of partitions $\mathcal{P}$ stabilizes.}{  
				$\mathcal{P}  \leftarrow \mathcal{P} \bigvee \partition \left(f_{\mathcal{L}}( \mathcal{B}_\mathcal{P}) \right)$   \\
				$\mathcal{P}  \leftarrow \mathcal{P} \bigvee  \partition \left(f_{X^2}( \mathcal{B}_\mathcal{P})\right) $ 
			}
		\end{algorithm}
	\end{subfigure}\\
        \caption{Algorithms for finding bases $\mathcal{B}_\mathcal{R} \subset
        \mathbb{S}^n$ and $\mathcal{B}_\mathcal{P} \subset \mathbb{S}^n$ of
        $\mathcal{S}_{coord}$ and $\mathcal{S}_{part}$, respectively. One
        algorithm  grows a relation $\mathcal{R} \subseteq [n] \times [n]$ and
        the other  refines a partition $\mathcal{P}$ of $[n] \times [n]$. (Here, $\mathcal{P}_1 \bigvee \mathcal{P}_2$ denotes the coarsest common refinement of partitions $\mathcal{P}_1$ and $\mathcal{P}_2$.) The
        inputs are $C_{\mathcal{L}}, Y_{\mathcal{L}^{\perp}} \in \mathbb{S}^n$
        and the linear subspace $\mathcal{L} \subseteq \mathbb{S}^n$.
        }\label{fig:algComb}
\end{figure}

\subsection{Randomization via sampling}\label{sec:partrand}

Explicitly constructing symbolic representations of $f_{\mathcal{L}}( \mathcal{B})$ and $f_{X^2}( \mathcal{B})$
is not necessary for finding the partitions and relations they induce. One can instead evaluate the maps $X \mapsto X^2$ and $P_{\mathcal{L}}:\mathbb{S}^n\rightarrow \mathbb{S}^n$ at a random combination of  elements in $\mathcal{B}$. Consider, for instance,
a point evaluation of $f_{X^2}( \mathcal{B}_\mathcal{P})$ at $t^* \in \mathbb{R}^{|\charmats|}$, i.e., consider
\[
  {f_{X^2}(\charmats)}_{|t=t^*} :=  {\left(\sum_{B \in \charmats } t^*_{B} B \right)}^2.
\]
The supports of $f_{X^2}(\charmats)$ and ${f_{X^2}(\charmats)}_{|t=t^*}$ are the same
for almost all  $t^*$. Similarly, the partitions induced by
$f_{X^2}(\charmats)$ and ${f_{X^2}(\charmats)}_{|t=t^*}$ are the same, i.e.,  
\[
  \partition \left(f_{X^2}(\charmats)\right)=\partition \left({f_{X^2}(\charmats)}_{|t=t^*}\right),
\]
for almost all  $t^*$.  The following  illustrates this equality for a particular $t^*$.

\addtocounter{ex}{-1}
\begin{ex}[continued]
  For $\mathcal{B}$  defined   previously,  the point evaluation ${f_{X^2}(\mathcal{B})}_{|t =t^*_{\mathcal{B}}}$ at $t^*_{\mathcal{B}} = (2,3,4)$ is
	\[
          {f_{X^2}(\mathcal{B})}_{|t =t^*_{\mathcal{B}}} =  {( 2U + 3 V + 4 W)}^2  = \left(\begin{array}{cccc} 20 & 0 & 8 & 6\\ 0 & 20 & 6 & 8\\ 8 & 6 & 13 & 0\\ 6 & 8 & 0 & 13 \end{array}\right).
	\]
        We see the partition $\partition \left({f_{X^2}(\mathcal{B})}_{|t =t^*_{\mathcal{B}}}\right)$ is the same as the partition $\partition \left(f_{X^2}(\charmats)\right)$ given by~\eqref{ex:Partition}.
\end{ex}

\subsection{Generalization of partition and coordinate subspaces}\label{sec:alg01}

Coordinate and partition subspaces have a trivial common generalization: subspaces
with an orthogonal basis of $0/1$ matrices, or, equivalently, a basis
of $0/1$ matrices with disjoint support. This motivates the following definition:
\[
\mathcal{S}_{0/1} :=  \bigcap \{  \mathcal{S} \subseteq \mathbb{S}^n : \mathcal{S} \mbox{ is admissible  and has an orthogonal basis of 0/1 matrices} \}.
\]
A procedure for finding $\mathcal{S}_{0/1}$ combines features of the algorithms presented for finding $\mathcal{S}_{part}$ and $\mathcal{S}_{coord}$: specifically, it iteratively grows a relation $\mathcal{R}$ and refines a partition  $\mathcal{P}$ of this relation.
	\begin{center}  \qquad \qquad \qquad \qquad
		\begin{minipage}{.65 \textwidth} \DontPrintSemicolon
			Initialize $\mathcal{R}$ to $\relation (C_{\mathcal{L}})  \bigcup  \relation (Y_{\mathcal{L}^{\perp}})$ \\
			Initialize $\mathcal{P}$ to $\partition_{\mathcal{R}} (C_{\mathcal{L}}) \bigvee \partition_{\mathcal{R}} (Y_{ \mathcal{L}^{\perp}})$
			\\
                        \Repeat{\rm ascending chain of subspaces $\Span (\mathcal{B}_{\mathcal{P}})$ stabilizes.}{
				\For{$f \in \{f_{\mathcal{L}}, f_{X^2}\}$}{
					Replace $\mathcal{R}$ with  $\mathcal{R} \cup \relation \left(f( \mathcal{B}_\mathcal{P}) \right)$    \\
					Add class $\mathcal{R} \setminus \left( \cup_{P \in \mathcal{P} } P \right)$ to  $\mathcal{P}$   \\
					Replace $\mathcal{P}$ with refinement   $\mathcal{P} \bigvee \partition_{\mathcal{R}} \left( f( \mathcal{B}_\mathcal{P})\right)$ 
				}
			}  
		\end{minipage}
	\end{center} 
Here $\partition_{\mathcal{R}}(T)$ denotes the partition of $\mathcal{R} \subseteq [n] \times [n]$ induced by the unique entries of a matrix $T$  with support contained in $\mathcal{R}$.

\section{Examples}\label{sec:example}

We now apply our techniques to SDPs rising in applications.  We do all
computation on an Intel 3GHz desktop with 128~gigabytes of RAM.\@ The
algorithms of Theorem~\ref{thm:subalg}~and~Figure~\ref{fig:algComb}, used to
identify an admissible subspace $\mathcal{S}$, and the algorithm of~\cite[Chapter
6]{permenterThesis2017}, used to find the linear map $\isoName$ and cone
$\coneNameTwo$ satisfying $\isoname(\coneNameTwo) = \mathbb{S}^n_{+} \cap
\mathcal{S}$, were all implemented in MATLAB.\@  We use the solver SeDuMI~\cite{sturm1999using}
to solve the SDPs.

\paragraph{Format of original SDPs}
Each primal-dual pair is originally expressed in either
SeDuMi~\cite{sturm1999using} or SDPA~\cite{fujisawa2002sdpa} format and may
have a mix of free and conic variables, where the cones are either nonnegative
orthants or cones of psd matrices.\footnote{These formats also allow for
Lorentz cones. None of the examples presented, however, use this type of cone.}
From these formats, we  eliminate free variables, reformatting the primal
problem as
\begin{align}\label{sdp:mainEq}
\begin{array}{ll}
\mbox{ minimize } & \langle C, X \rangle \\
\mbox{ subject to } & \langle A_i, X \rangle = b_i \;\;\; \forall i \in [m]   \\
& X \in \mathbb{S}_{+}^{n_{1}} \times \cdots \times \mathbb{S}_{+}^{n_{r}},   
\end{array} 
\end{align}
where $C, A_i \in  \mathbb{S}^{n_{1}} \times \cdots \times \mathbb{S}^{n_{r}}$
are fixed, and $\langle \cdot, \cdot \rangle$ denotes the inner-product
obtained by equipping each $\mathbb{S}^{n_{i}}$ with the trace
inner-product. (This reformatting amounts to eliminating free variables and relabeling nonnegative orthants 
as  products of psd cones of order one.)
We will in some cases report the number of non-zero (nnz) entries in a
description of~\eqref{sdp:mainEq}; this equals the number of non-zero
floating-point numbers needed to store $C$ and ${\{A_i\}}^m_{i=1}$.  We also
report a tuple of  ranks for~\eqref{sdp:mainEq}, which is simply the tuple
$(n_1,\ldots,n_r)$.

\paragraph{Format of reformulations} 
\renewcommand{\coneName}{\mathcal{C}}
We reformulate each SDP by finding an admissible subspace $\mathcal{S}$, simple algebras $\mathcal{J}_i$,
and an injective homomorphism $\isoName : \oplus^q_{i=1} \mathcal{J}_i \rightarrow \mathbb{S}^n$ satisfying
\[
\mathbb{S}^n_{+}  \cap\mathcal{S}  = \isoName(\coneName_1 \times \cdots \times \coneName_q),
\]
where $\coneName_i$ is the cone-of-squares of $\mathcal{J}_i$. The reformulation is as in
Proposition~\ref{prop:equiv2}:
\begin{align}\label{sdp:optionTwo}
\begin{array}{ll}
\mbox{ minimize } & \langle \isoName^*(C), \hat X \rangle \\
\mbox{ subject to } & \langle \isoName^*(A_i), \hat X \rangle =  b_i \;\;\; \forall i \in T \subseteq [m]   \\
& \hat X \in \coneName_1 \times \cdots \times \coneName_q,
\end{array} 
\end{align}
where $T$ indexes a maximal linearly-independent subset of equations.
We will in some cases report the number of non-zero (nnz) entries in a
description of~\eqref{sdp:optionTwo}; this equals the number of non-zero
floating-point numbers needed to store $\isoName^*(C)$ and ${\{\isoName^*(A_i)
\}}_{i \in T}$.  We also report the tuple $(r_1,\ldots,r_q)$, where $r_i$ is the rank of $\mathcal{J}_i$.

\begin{rem}
        For most examples, $\coneName_1 \times \cdots \times \coneName_q$ is  a
        product of psd cones  $\mathbb{S}^{r_1}_{+} \times \cdots \times
        \mathbb{S}^{r_q}_{+}$  and the tuple $(r_1,\ldots,r_q)$ indicates their
        orders.   We discuss the
        only exception in Section~\ref{sec:acomplexalg}.  We also note
        $\mathbb{S}^2$ is isomorphic to a spin-factor algebra---hence,
        $\mathbb{S}_{+}^2$ is isomorphic to a Lorentz cone.
\end{rem}

\paragraph{Reference subspaces and inclusions}
For convenience, we will let $\mathcal{S}_{full} := \mathbb{S}^{n_{1}} \times
\cdots \times \mathbb{S}^{n_{r}}$ denote the full ambient space of the original
instance~\eqref{sdp:mainEq}. As discussed in~\cite{de2010exploiting}, an SDP
can be restricted to the *-algebra generated by its data matrices; see also Proposition~\ref{prop:staralgdata}.  To compare
with this restriction,  we   let
\[
\mathcal{S}_{data}:=\mathcal{M}_{data} \cap (\mathbb{S}^{n_{1}} \times \cdots \times \mathbb{S}^{n_{r}}),
\] 
where $\mathcal{M}_{data}\subseteq \mathbb{R}^{n \times n}$ is the  *-subalgebra of $\mathbb{R}^{n \times n}$ generated by the problem data $C$ and ${\{A_i\}}^m_{i=1}$ (using matrix multiplication as a product and transposition as the *-involution). Recall that 
\[
\mathcal{S}_{opt} \subseteq  \mathcal{S}_{0/1} \subseteq  \mathcal{S}_{coord} \subseteq \mathcal{S}_{full}.
\]
We'll see that different inclusions  hold strictly for different examples. By definition, we also have that 
\[
\mathcal{S}_{opt} \subseteq \mathcal{S}_{data} \subseteq \mathcal{S}_{full}.
\]
Examples will show that $\mathcal{S}_{opt}$ can be (much) smaller than  $\mathcal{S}_{data}$.

\subsection{Libraries of problem  instances}

The first set of SDPs come from three public sources:  the
parser SOSTOOLS~\cite{papachristodoulou2013sostools}, the DIMACS
library~\cite{pataki1999dimacs} and a set of structured SDP instances
from~\cite{de2009new}.    Table~\ref{tab:library} reports the dimensions of the
subspaces $\mathcal{S}_{opt}$, $\mathcal{S}_{0/1}$,  $\mathcal{S}_{coord}$,
$\mathcal{S}_{data}$ and $\mathcal{S}_{full}$.  Note   the  inclusions
$\mathcal{S}_{opt} \subseteq  \mathcal{S}_{0/1} \subseteq \mathcal{S}_{coord}
\subseteq \mathcal{S}_{full}$  hold as expected, and, as
Table~\ref{tab:library} indicates, different ones hold strictly for different
instances. For a large fraction, $\mathcal{S}_{full}$ equals
$\mathcal{S}_{data}$, implying  generating a *-subalgebra from the problem
data~\cite{de2010exploiting}  does not simplify these instances.

\begin{rem} 
	We note the  libraries~\cite{pataki1999dimacs,de2009new} have additional instances on which
	our method was not effective $(\mathcal{S}_{opt} = \mathcal{S}_{full}$); we do not
	report results for these instances.  
\end{rem}

\begin{rem}
The kissing number and copositivity instances of Table~\ref{tab:library}
  (denoted {\tt kissing\_x\_y\_z} and {\tt coposxy}) can also be simplified using
  group theoretic techniques 
  \cite{caluza2018improving,
  dobre2015exploiting} (related to Proposition~\ref{prop:csigroup}) that are tailored to these specific SDP families.
\end{rem}

\begin{table}  
	\begin{center}
		\begin{tabular}{ccccccc}  
			instance & $\mathcal{S}_{opt}$ & $\mathcal{S}_{0/1}$ & $\mathcal{S}_{coord}$ &  $\mathcal{S}_{data}$ &  $\mathcal{S}_{full}$ & References\\   \toprule
			{\tt sosdemo2}  &  25 & 25 & 28 & 103 & 103  & \multirow{7}{*}{\parbox{1.5cm}{Instances from~\cite{papachristodoulou2013sostools}   }  }  \\ 
			{\tt sosdemo4} &  11 & 11 & 85 & 630 & 630 & \\ 
			{\tt sosdemo5}  & 226 &816 & 816 & 816 & 816 & \\ 
			{\tt sosdemo6} &  49 & 49 & 327 & 462 & 462 & \\ 
			{\tt sosdemo7} &  40 & 40 & 68 & 68 & 68 \\ 
			{\tt sosdemo9} &  26 & 26 & 26 & 78 & 78 &\\ 
			{\tt sosdemo10} &  78 & 78 & 78 & 254 & 254&  \\  \midrule
			{\tt hamming\_7\_5\_6} &  5 & 5 & 8256 & 8256 & 8256  &  \multirow{8}{*}{\parbox{1.5cm}{Instances from~\cite{pataki1999dimacs}  }  } \\ 
			{\tt hamming\_8\_3\_4} &  5 & 5 & 32896 & 32896 & 32896 & \\
			{\tt hamming\_9\_5\_6} &  6 & 6 & 131328 & 131328 & 131328 & \\
			{\tt hamming\_9\_8}  &  6 & 6 & 131328 & 131328 & 131328 & \\ 
			{\tt hamming\_10\_2}   &  7 & 7 & 524800 & 524800 & 524800 &\\ 
			%{\tt hamming\_11\_2}   &  7 & 7 & 2098176 & 2098176 & 2098176 &\\ 
			{\tt copo14}  &  73 & 73 & 1834 & 1834 & 1834 &  \\ 
			{\tt copo23} &  188 & 188 & 8119 & 8119 & 8119 & \\ 
			{\tt copos68}   & 1576   & 1576  & 209644  &  209644 &209644  &  \\
			\midrule 
			{\tt ThetaPrimeER23\_red}  &  86 & 762 & 777 & 101 & 1712 &  \multirow{7}{*}{\parbox{1.5cm}{Instances from~\cite{de2009new}  }  } \\ 
			{\tt ThetaPrimeER29\_red}  &  104 & 1125 & 1143 & 122 & 2486 & \\ 
			{\tt ThetaPrimeER31\_red}  &  110 & 1262 & 1281 & 129 & 2776 & \\ 
			{\tt crossing\_K\_7n} &  113 & 577 & 3138 & 113 & 3138 & \\ 
			{\tt crossing\_K\_8n} &  479 & 18577 & 72630 & 479 & 72630 & \\ 
			{\tt kissing\_3\_5\_5} &  811 & 811 & 3796 & 3796 & 3796 &  \\ 
			{\tt kissing\_4\_7\_7} &  3723 & 3723 & 19760 & 19760 & 19760 &
		\end{tabular} 
		\caption{Dimensions of admissible subspaces $\mathcal{S}_{opt}$, $\mathcal{S}_{0/1}$
			and $\mathcal{S}_{coord}$ compared with dimensions of the ambient space $\mathcal{S}_{full}$ and $\mathcal{S}_{data}$---the 
                        (symmetric part) of the *-algebra generated by   $C$ and ${\{A_i\}}^m_{i=1}$.}\label{tab:library}
		\vspace{1cm}
		\begin{tabular}{ccccc} 
                  instance & $\mathcal{S}_{opt}$ & $\mathcal{S}_{0/1}$ & $\mathcal{S}_{coord}$ &  $\mathcal{S}_{full}$ \\ \toprule
			{\tt ThetaPrimeER23\_red} & ($3,2_{12\times},1_{44\times})$ & ($27,25,5,1_{44\times})$ & ($27,25,5,1_{59\times})$  & ($57,1_{59\times})$ \\ 
			{\tt ThetaPrimeER29\_red} & ($3,2_{15\times},1_{53\times})$ & ($33,31,5,1_{53\times})$ & ($33,31,5,1_{71\times})$  & ($69,1_{71\times})$ \\ 
			{\tt ThetaPrimeER31\_red} & ($3,2_{16\times},1_{56\times})$ & ($35,33,5,1_{56\times})$ & ($35,33,5,1_{75\times})$  & ($73,1_{75\times})$   \\ 
                        {\tt crossing\_K\_7n} & $(3_{6 \times}, 2_{4 \times}, 1_{65 \times})$ & ($13_{4 \times},1_{57\times}$) & ($79,1_{57\times}$)  & ($79,1_{57\times}$)  \\ 
			{\tt crossing\_K\_8n} &$(7_{2\times},5_{2\times},4_{9\times},3_{7\times},2_{4\times},1_{249\times})$ &$(105,97,92,86,1_{240\times})$ & $(380,1_{240\times})$ & $(380,1_{240\times})$  \\ 
		\end{tabular}  
		
                \caption{Tuple of ranks  for select examples after restricting to indicated subspace. Here,  $s_{t \times}$ means $s$  repeated $t$ times, i.e., $3_{2\times} := (3,3)$.     }\label{tab:libraryinstances}

                \vspace{1cm} 
		\begin{tabular}{ccccc} 
                  instance & $t_{orig}$ & $\mathcal{S}_{opt}$ & $\mathcal{S}_{0/1}$ & $\mathcal{S}_{coord}$  \\ \toprule
                        {\tt ThetaPrimeER23\_red} & 0.21 & 0.16, 0.09  & 0.12, 0.12  & 0.02, 0.13 \\ 
                        {\tt ThetaPrimeER29\_red} & 0.19 & 0.14, 0.10  & 0.11, 0.15  & 0.02, 0.14 \\ 
                        {\tt ThetaPrimeER31\_red} & 0.25 & 0.17, 0.15  & 0.13, 0.19  & 0.02, 0.19 \\ 
                        {\tt crossing\_K\_7n}     & 0.33 & 0.25, 0.12  & 0.18, 0.14  & 0.02, 0.33 \\ 
                        {\tt crossing\_K\_8n} & 56.7     & 2.48, 0.58  & 2.34, 10.37  & 0.02, 56.7 \\
		\end{tabular} \;\;
		\begin{tabular}{ccc} 
                  instance & $t_{orig}$  & $\mathcal{S}_{opt}$  \\ \toprule
                  {\tt hamming\_7\_5\_6} & 10.12   &  0.09, 0.04 \\   
                  {\tt hamming\_8\_3\_4} & 4 hours &  0.15, 0.02 \\
                  {\tt hamming\_9\_5\_6} & Fail     &  0.48, 0.02 \\
                  {\tt hamming\_9\_8\  } & Fail     &  0.49, 0.02 \\
                  {\tt hamming\_10\_2 }  & Fail     &  2.29, 0.03 \\
                  %{\tt hamming\_11\_2 }  & Fail     &  25.34, 0.08 \\
		\end{tabular} \ 
                \caption{The original solver time $t_{orig}$ (in seconds) and a list $t_{pre}, t_{solve}$
                of preprocessing and solver times for restrictions to indicated subspaces. Failures were due to
                insufficient memory.}\label{tab:hammingSolve}
        \end{center}
\end{table}

\subsubsection{The {L}ov{\'a}sz number} 
We give special attention to the  Table~\ref{tab:libraryinstances} instances  
denoted ${\tt hamming\_m\_x}$ and ${\tt hamming\_m\_x\_y}$, taken from~\cite{pataki1999dimacs}.
For a specific graph $G$ with vertices $\{1,\ldots,n\}$ and edge set $E$, each instance has the following form
\begin{align}\label{sdp:theta}
&\begin{array}{lll}
\mbox{ maximize} &   \trace \ones  X   \\ 
\mbox{ subject to}  & \trace X  = 1 , X \in \mathbb{S}^n_{+} \\
& \trace (E_{ij} + E_{ji}) X = 0 \qquad \forall (i,j) \in E, 
\end{array} 
\end{align}
where $\ones \in \mathbb{S}^{n}$ is the all-ones-matrix and $E_{ij}$ is a standard basis matrix of $\mathbb{R}^{n \times n}$.    

The  graphs for these instances are closely related to the Hamming graph
$H(m,d)$, whose nodes are the Boolean vectors of length $m$ that  are adjacent
iff their  Hamming distance is at least $d$.  The graphs of  ${\tt
hamming\_m\_x}$ and ${\tt hamming\_m\_x\_y}$  are modifications of such graphs:  nodes   are adjacent iff their  Hamming distance  is {\tt
x} or is {\tt x} or {\tt y}. When $G$ is a Hamming graph, it is well known that
one can convert SDP~\eqref{sdp:theta}  into a linear program   using the theory
of association schemes~\cite{schrijver1979comparison}. Unsurprisingly, we find similar
simplifications for the modified graphs; precisely, $\mathbb{S}^n_{+}
\cap \mathcal{S}_{opt}$ is isomorphic to a non-negative orthant of order equal
to the dimension of $\mathcal{S}_{opt}$, i.e.,  
\[
\mathbb{S}_{+}^n \cap \mathcal{S}_{opt} = \isoName( \mathbb{R}^{\dim \mathcal{S}_{opt} }_{+})
\]
for an injective map $\isoName : \mathbb{R}^{\dim \mathcal{S}_{opt}} \rightarrow \mathbb{S}^n$.  

As reported~\cite{mittelmann2003independent}, these instances are challenging
for a wide array of solvers due to their size; indeed, we are only able to
solve two of them directly (Table~\ref{tab:hammingSolve}).  Constructing the
reformulation over $\mathcal{S}_{opt}$, however, converts each SDP into a
trivial linear program. Further, finding $\mathcal{S}_{opt}$ and constructing
the reformulation takes negligible effort compared to original solver time
(Table~\ref{tab:hammingSolve}).  Note the other automated approach---generating
a $*$-algebra from the data matrices $\ones$, $I$, and
${\{E_{ij} + E_{ji}\}}_{ (i,j) \in E }$---fails for these instances, i.e., $\mathcal{S}_{data} = \mathcal{S}_{full}$ (Table~\ref{tab:libraryinstances}).

\subsubsection{Decompositions and majorization}
In Table~\ref{tab:libraryinstances} we report the tuple of ranks for
the subspaces $\mathcal{S}_{opt}$, $\mathcal{S}_{0/1}$  and $\mathcal{S}_{coord}$
for select examples to confirm our main theorem on optimal decompositions (Theorem~\ref{thm:OptimalDecomp}).
Specifically, we select examples satisfying the strict inclusions:
\[
\mathcal{S}_{opt} \subset  \mathcal{S}_{0/1} \subset  \mathcal{S}_{coord}.
\]
Given these strict inclusions, Theorem~\ref{thm:OptimalDecomp} predicts  the
ranks of $\mathcal{S}_{0/1}$ and  $\mathcal{S}_{coord}$  weakly majorize  those
of $\mathcal{S}_{opt}$ in the sense of Definition~\ref{defn:majorize}.
Similarly, it predicts  the ranks of $\mathcal{S}_{coord}$  weakly majorize
those  of $\mathcal{S}_{0/1}$.

Table~\ref{tab:libraryinstances} confirms both these predictions. The first row, for instance, reports the
following tuples  $r_1  \in \mathbb{Z}^{l_1}$ and $r_2 \in \mathbb{Z}^{l_2}$
for $\mathcal{S}_{opt}$ and $\mathcal{S}_{0/1}$, respectively:  
\[
r_1 := (3,\underbrace{2,2,\ldots,2}_{12 \times},\underbrace{1,1,\ldots,1}_{44 \times})  \qquad r_{2}  :=(27,25,5,\underbrace{1,1,\ldots,1}_{44\times}).
\]
It easily follows  $r_{2}$ weakly majorizes $r_1$, i.e., for all positive integers $q \in \mathbb{Z}$,  
\[
  \sum^{\min \{q,l_{2 }  \}}_{i=1} {[r_{2}]}_i \ge \sum^{\min \{q,l_{1 }  \}}_{i=1} {[r_{1}]}_i.
\]

As also expected, for the instances of Table~\ref{tab:libraryinstances},
reformulating over $\mathcal{S}_{opt}$ reduces solver time the most, but
requires the most preprocessing (Table~\ref{tab:hammingSolve}).  In fact, for
some instances, solver time reductions do not offset the extra preprocessing time, justifying the larger (but easier
to construct) reformulations over $\mathcal{S}_{0/1}$ and
$\mathcal{S}_{coord}$. (To further reduce preprocessing time, Section~\ref{ex:cprank} introduces an alternative reformulation~\eqref{sdp:optionOne} to~\eqref{sdp:optionTwo}
that reduces the dimension of the feasible set but doesn't simplify the cone constraint.)

\subsubsection{An algebra with a complex direct-summand}\label{sec:acomplexalg}

The example ${\tt sosdemo5}$ is an SDP that bounds a quantity from robust control theory---the structured singular value $\mu(M,\mathbf{\Delta})$~\cite{packard1993complex}:
\begin{align}
\mu(M,\mathbf{\Delta}) := \frac{1}{\inf \{ \| \Delta \| : \Delta \in \mathbf{\Delta}, \det(I-M\Delta) = 0   \} }.
\end{align}
Here, $M$ is a complex matrix and $\mathbf{\Delta}$ is a set of complex
matrices. Though the parameters of $\mu(M,\mathbf{\Delta})$ are complex, one
can formulate an SDP with real data matrices to bound $\mu(M,\mathbf{\Delta})$.
This is done in ${\tt sosdemo5}$  for particular $M$ and $\mathbf{\Delta}$.
After decomposing $\mathcal{S}_{opt}$ into a direct-sum of minimal
ideals, we find one of the direct-summands is isomorphic to an algebra of
complex Hermitian  matrices. Precisely, $\mathcal{S}_{opt} = \oplus^{11}_{i=1}
\mathcal{S}_i$ for minimal ideals $\mathcal{S}_i$. Letting $r:=(\rank
\mathcal{S}_{1},\ldots,\rank \mathcal{S}_{11})$ and $d:=(\dim \mathcal{S}_{1},\ldots,\dim \mathcal{S}_{11})$, we have
\begin{align*}
r &= (1  ,   1 ,    1 ,    1 ,    4  ,   4^* ,    4 ,    6  ,  10 ,   10    ,10) \\
d &= (1  ,   1 ,    1  ,   1 ,   10   , 16^*   , 10  ,  21 ,   55 ,   55   , 55).
\end{align*}
Note with the exception of the entries marked $*$, the relation $d_i = { r_i+1 \choose 2}$
holds, showing $\mathcal{S}_i$ is isomorphic to the algebra of real symmetric matrices of order $r_i$.
The exception satisfies $d_i = r^2_i$, showing the corresponding ideal $\mathcal{S}_i$  is isomorphic to the algebra of complex Hermitian matrices of order $r_i$. We remark this is the only example considered where the direct-summands are not all isomorphic to $\mathbb{S}^n$ for some $n$.

\subsection{Coordinate subspaces and sparse decompositions}\label{ex:sosopt}

We next consider SDPs  constructed by demonstration scripts packaged with the control system analysis tools available at  
\[
{\tt http://www.aem.umn.edu/\textasciitilde AerospaceControl/ },
\]
which build upon   the parser SOSOPT~\cite{seiler2013sosopt}.   For these SDPs, the optimal subspace
$\mathcal{S}_{opt}$ equals the optimal coordinate subspace
$\mathcal{S}_{coord}$. As indicated in Table~\ref{tab:sosopt}, these SDPs illustrate   we can always restrict to
$\mathcal{S}_{coord}$ without increasing the number of non-zero entries in the
problem description, since restricting to $\mathcal{S}_{coord}$ amounts to setting
certain off-diagonal entries of the data to zero. Though these examples are of
small size, they illustrate  $\mathcal{S}_{coord}$ is a proper subspace of
$\mathcal{S}_{full}$ for many SDPs arising in sums-of-squares
optimization.

\begin{rem}
	Note some of these scripts construct more than one SDP;\@ reported  results  are for the \emph{first} SDP constructed.
\end{rem}

\begin{table} \centering
	\begin{tabular}{ccccc}
		& \multicolumn{2}{c}{Orig.}  &  \multicolumn{2}{c}{$\mathcal{S}_{coord}$ }      \\ 
		&   ranks &  nnz & ranks &  nnz        \\
		\toprule 
                {\tt Chesi[1|4]\_IterationWithVlin} & $(9,5)$ & 181 & $(6,3_{2\times},2)$ & 97 \\ 
		{\tt Chesi3\_GlobalStability} & $(14,5)$ & 341 & $(8,6,3,2)$ & 193 \\ 
                {\tt Chesi[5|6]\_Bootstrap} & $(19,9)$ & 928 & $(13,6_{2\times},3)$ & 520 \\  
		{\tt Chesi[5|6]\_IterationWithVlin} & $(19,9)$ & 928 & $(13,6_{2\times},3)$ & 520 \\ 
		{\tt Coutinho3\_IterationWithVlin} & $(9,5)$ & 181 & $(6,3_{2\times},2)$ & 97 \\ 
		{\tt HachichoTibken\_Bootstrap} & $(19,9)$ & 685 & $(12,7,6,3)$ & 373 \\ 
		{\tt HachichoTibken\_IterationWithVlin} & $(19,9)$ & 685 & $(12,7,6,3)$ & 373 \\ 
		{\tt Hahn\_IterationWithVlin} & $(9,5)$ & 156 & $(6,3_{2\times},2)$ & 84 \\ 
		{\tt KuChen\_IterationWithVlin} & $(19,9)$ & 928 & $(13,6_{2\times},3)$ & 520 \\ 
		{\tt  Parrilo1\_GlobalStabilityWithVec} & $(3,2)$ & 20 & $(2,1_{3\times})$ & 14 \\ 
		{\tt  Parrilo2\_GlobalStabilityWithMat} & $(3,2)$ & 16 & $(2,1_{3\times})$ & 10 \\ 
		{\tt Pendubot\_IterationWithVlin} & $(14,4)$ & 372 & $(10,4_{2\times})$ & 292 \\ 
		{\tt  VDP\_IterationWithVball} & $(5,4)$ & 82 & $(3_{2\times},2,1)$ & 55 \\ 
		{\tt  VDP\_IterationWithVlin} & $(9,5)$ & 181 & $(6,3_{2\times},2)$ & 97 \\ 
		{\tt  VDP\_LinearizedLyap} & $(9,5)$ & 156 & $(6,3_{2\times},2)$ & 84 \\ 
		{\tt  VDP\_MultiplierExample} & $(5,2)$ & 37 & $(3,2,1_{2\times})$ & 23 \\ 
		{\tt  VannelliVidyasagar2\_Bootstrap} & $(19,9)$ & 928 & $(13,6_{2\times},3)$ & 520 \\ 
		{\tt  VannelliVidyasagar2\_IterationWithVlin} & $(19,9)$ & 928 & $(13,6_{2\times},3)$ & 520 \\ 
		{\tt  VincentGrantham\_IterationWithVlin} & $(9,5)$ & 181 & $(6,3_{2\times},2)$ & 97 \\ 
		{\tt WTBenchmark\_IterationWithVlin} & $(19,9)$ & 685 & $(13,6_{2\times},3)$ & 385 \\ 
		
		%copos\_1 & $35$ & 1225 & $(5_{5\times},1_{10\times})$ & 135 \\ 
		%copos\_2 & $120$ & 14400 & $(8_{8\times},1_{56\times})$ & 568 \\ 
		%copos\_3 & $286$ & 81796 & $(11_{11\times},1_{165\times})$ & 1496 \\ 
		%copos\_4 & $560$ & 313600 & $(14_{14\times},1_{364\times})$ & 3108 \\ 
		%copos\_5 & $969$ & 938961 & $(17_{17\times},1_{680\times})$ & 5593 \\ 
        \end{tabular}   \caption{Ranks  and number of non-zero (nnz) entries in problem description  of original instance  and its  restriction~\eqref{sdp:optionTwo} to  $\mathcal{S}_{coord}$. The notation $r_{s \times}$ indicates $r$ repeated $s$ times. 
           }\label{tab:sosopt}
\end{table}

\subsection{Comparison with LP method of~\citeauthor*{grohe2014dimension} }
In~\cite{grohe2014dimension},~\citeauthor{grohe2014dimension} describe a
reduction method for \emph{linear} programming (LP) and show it outperforms a
symmetry reduction method of~\cite{bodi2011symmetries} on a collection of LPs;
indeed, they show their method theoretically
subsumes~\cite{bodi2011symmetries}.  The linear programs used for comparison
are  relaxations of    integer programs  studied
in~\cite{margot2003exploiting}. By treating each linear inequality as a
semidefinite constraint of order one, we  applied our method to the same LP
relaxations. Of the 57 relaxations, we find the same reductions on 56.  For the
remaining instance ({\tt  cov1054sb}), we 
outperform~\cite{grohe2014dimension}.  For space reasons, Table~\ref{tab:lp}
reports results for just a small subset of these LP relaxations.  To
match~\cite{grohe2014dimension}, we give the number of \emph{dual} variables
and  inequality constraints.  In terms of     SDP~\eqref{sdp:mainEq} and the
SDP~\eqref{sdp:optionTwo}, the number of dual variables and constraints equals
the number of linear equations and the sum of the ranks, respectively.

%That reductions are the same across the other examples is interesting  since the theoretical
%development of  our method and theirs  is different, at least superficially. 

\begin{table} \centering
	\begin{tabular}{ccccccc} & \multicolumn{3}{c}{Constraints} & \multicolumn{3}{c}{Variables}  \\ 
		& Orig. & CR  & $\mathcal{S}_{opt}$ & Orig. & CR  & $\mathcal{S}_{opt}$ \\ \toprule
		{\tt cov1053 } & 252 & 1 & 1 & 679 & 5  & 5    \\
		{\tt cov1054} & 252 & 1 & 1 & 889 & 6  & 6     \\
		{\tt cov1054sb}& 252 & 252 & 1 & 898 & 898  &  6      \\
		{\tt cov1075}  & 120 & 1 & 1 & 877 & 7  & 7         \\ 
		{\tt cov1076}  & 120 & 1 & 1 & 835 & 7  & 7     \\ 
		{\tt cov1174}  & 330 & 1 & 1 & 1221 & 6  & 6    \\ 
		{\tt cov954}   & 126   & 1 & 1 & 507 & 6  & 6       
	\end{tabular}
	\caption{Dual variables and constraints of original LP, the LP  formulated via the color refinement (CR) method of~\cite{grohe2014dimension}, and the LP formulated via restriction to  $\mathcal{S}_{opt}$. Columns labeled (CR) use numbers reported in~\cite{grohe2014dimension}.}\label{tab:lp}
\end{table}

\subsection{Completely-positive rank, the subspace $\mathcal{S}_{0/1}$, and decomposition trade-offs}\label{ex:cprank}
Our last example  illustrates   restrictions to $\mathcal{S}_{0/1}$,
the optimal subspace with an  orthogonal basis of $0/1$ matrices.   The considered SDP family    yields lower-bounds  of \emph{completely-positive  rank}, or cp-rank for short.  The cp-rank of $W \in \mathbb{S}^n_{+}$ measures the size of the smallest non-negative factorization of $W$.  Precisely, it is the smallest 
$r$ for which $V \in \mathbb{R}_{+}^{n \times r}$  exists satisfying  $W = VV^T$. (It is infinite if such a factorization does not exist for any $r$.)
As shown in~\cite{fawzi2014self}, the cp-rank of $W \in \mathbb{S}^n$ is lower bounded by the optimal value of the following SDP:\@
\begin{align*}
\begin{array}{crl}  \mbox{minimize } t \\ \mbox{subject to} \\ & \left( \begin{array}{cc} t & \vect W^T \\ \vect W & X \end{array} \right) \succeq 0 \\ 
& X_{ij,ij} \le W^2_{ij} & \forall i,j \in \{1,\ldots,n\}  \\
& X \preceq W \otimes W \\
& X_{ij,kl} = X_{il,jk} & \forall (1,1) \le (i,j) < (k,l) \le (n,n).
\end{array} 
\end{align*}
Here, $W \otimes W$ denotes the \emph{Kronecker product} and $\vect W$ denotes the $n^2 \times 1$ vector obtained by stacking the columns of $W$. The double subscript $ij$ indexes the $n^2$ rows (or columns) of $X$ 
and the inequalities on $(i,j)$ and $(k,l)$  hold iff they hold element-wise. (See~\cite{fawzi2014self} for  clarification on this notation.)

In this example, we solve three instances of this SDP taking $W$ equal to the matrices $Z$,  $Z \otimes Z$,   and   $Z\otimes Z \otimes Z$,
where
\[
Z = \left(\begin{array}{ccc} 4 & 0 & 1 \\ 0 & 4 & 1 \\ 1 & 1 & 3 \end{array}\right).
\]
Table~\ref{cp:rank}  reports computational savings obtained by restricting to $\mathcal{S}_{0/1}$.

\paragraph{Alternative reformulation}
For   these examples, we compare~\eqref{sdp:optionTwo}
against an  alternative reformulation that reduces the dimension of the dual feasible set, but leaves the cone constraint unchanged.
It takes the following form
\begin{align}\label{sdp:optionOne}
\begin{array}{ll}
\mbox{ minimize } & \langle P_{\mathcal{S}_{0/1}}(C),   X \rangle \\
\mbox{ subject to } & \langle P_{\mathcal{S}_{0/1}}(A_i),   X \rangle =  b_i \;\;\; \forall i \in T \subseteq [m]   \\
&  X \in \mathbb{S}_{+}^{n_{1}} \times \cdots \times \mathbb{S}_{+}^{n_{r}},
\end{array} 
\end{align}
where $T \subseteq [m]$ indexes a maximal subset of linearly-independent  equations, and has dual  
\begin{align*}  
\begin{array}{ll}
\mbox{ maximize } & \sum_{i \in T } y_i b_i \\
\mbox{ subject to } & P_{\mathcal{S}_{0/1}}(C)- \sum_{i \in T } P_{\mathcal{S}_{0/1}}(A_i)   \in \mathbb{S}_{+}^{n_{1}} \times \cdots \times \mathbb{S}_{+}^{n_{r}}.
\end{array} 
\end{align*}
We can interpret the latter SDP as the dual of~\eqref{sdp:mainEq} restricted to the subspace $\mathcal{S}_{0/1}$,  recalling by Proposition~\ref{prop:optval}  that $\mathcal{S}_{0/1}$ contains both primal and dual solutions.  

Table~\ref{cp:rank}  shows solving~\eqref{sdp:optionOne} achieves computational savings and, indeed, can be preferred to solving~\eqref{sdp:optionTwo}.  As indicated, for the largest instance, we cannot even find the homomorphism $\isoname$ needed to construct~\eqref{sdp:optionTwo} due to memory constraints. For this example, the formulation~\eqref{sdp:optionOne} also preserves sparsity.

%      \begin{tabular}{   | c  | c | c  | c |c| }
%       \hline
%         Example &     Original     & $\mathcal{S}_{part}$ not decomposed     & $\mathcal{S}_{part}$ decomposed     \\
%      \hline
%   $Z$   &   $t_{solve} = 1.05$ & $t_{solve} = 0.16$, $t_{pre}=0.022$ &  $t_{solve} = 0.37$, $t_{pre}=0.30$ \\ \hline
%    $Z \otimes Z$  &  $t_{solve} = 39.17$ & $t_{solve} = 4.35$, $t_{pre}=0.5$ &  $t_{solve} = 2.09$, $t_{pre}=8.89$ \\  \hline
%    $Z \otimes Z \otimes Z$  & Out of memory  & $t_{solve} = 2008$, $t_{pre}=20$ &  Out of memory \\ \hline
%%        \end{tabular}
%    \caption{Solve time $t_{solve}$ for original SDP; Solve time $t_{solve}$ and   pre-processing time $t_{pre}$ for different versions of method.   }\label{tab:cprank3}
\begin{table}[tbp]
	\centering
	\subcaptionbox{Instance: $Z$}{%
		\begin{tabular}{cccccc}
			SDP   &    ranks & num eq &  $\nnz$ & $t_{pre}$ & $t_{solve}$  \\ \toprule
			Orig.~\eqref{sdp:mainEq}    & $(10,9,1_{9\times})$ & $37$ & $242$  & --- & .53   \\     
			Reform.~\eqref{sdp:optionTwo}     &   $(5,4,2_{4\times},1_{6\times})$ & $14$ & $859$  & 0.34 & 0.13  \\     
			Reform.~\eqref{sdp:optionOne}    &      $(10,9,1_{9\times})$ & $14$ & $242$ & 0.11 & 0.11  \\        
		\end{tabular}
	}  \qquad
	\subcaptionbox{Instance: $Z\otimes Z$}{%
		\begin{tabular}{cccccc}
			SDP      &    ranks &  num eq  &  $\nnz$ & $t_{pre}$ & $t_{solve}$ \\
			\toprule
			Orig.~\eqref{sdp:mainEq}     &   $(82,81,1_{81\times})$ & $2026$ & $15752$ & --- & 24.34  \\ 
			Reform.~\eqref{sdp:optionTwo}    &   $(12,11,10_{4\times},6_{4\times},4_{8\times},2_{2\times},1_{11\times})$ & $167$ & $158199$ & .98 & .90 \\    
			Reform.~\eqref{sdp:optionOne}       &  $(82,81,1_{81\times})$ & $167$ & $15752$  & 0.11 & 2.54
		\end{tabular}
	} \\
	\subcaptionbox{Instance: $Z\otimes Z \otimes Z$}{%
		\begin{tabular}{cccccc}
			SDP      &   ranks &  num eq  &  $\nnz$  &  $t_{pre}$ & $t_{solve}$ \\
			\toprule
			Orig.~\eqref{sdp:mainEq}    &  ($730,729,1_{729\times})$ & $142885$ & $1182290$  & \multicolumn{2}{c}{Out of memory}    \\   
			Reform.~\eqref{sdp:optionTwo}      &   \multicolumn{3}{c}{Out of memory} & \multicolumn{2}{c}{Out of memory}  \\          
			Reform.~\eqref{sdp:optionOne}     &  ($730,729,1_{729\times})$& $1883$ & $1182290$  & 6.5 & 1113    \\    
			
		\end{tabular}
	}    \caption{The first row  corresponds to the original SDP~\eqref{sdp:mainEq} and   the other rows to reformulations over $\mathcal{S}_{0/1}$.   Here, $t_{pre}$  is time spent (in seconds) finding  $\mathcal{S}_{0/1}$ and constructing the reformulation. Solve time $t_{solve}$ is also in seconds. 
}
\label{cp:rank}

\end{table}

\section*{Acknowledgements}

We thank Etienne de Klerk for useful discussions during the beginning stages of this  work. We also thank anonymous referees for comments that improved our presentation.

{\tiny
\bibliographystyle{abbrvnat}

\begin{thebibliography}{43}
\providecommand{\natexlab}[1]{#1}
\providecommand{\url}[1]{\texttt{#1}}
\expandafter\ifx\csname urlstyle\endcsname\relax
  \providecommand{\doi}[1]{doi: #1}\else
  \providecommand{\doi}{doi: \begingroup \urlstyle{rm}\Url}\fi

\bibitem[Alizadeh and Schmieta(2000)]{alizadeh2000symmetric}
F.~Alizadeh and S.~Schmieta.
\newblock Symmetric cones, potential reduction methods and word-by-word
  extensions.
\newblock In \emph{Handbook of Semidefinite Programming}, pages 195--233.
  Springer, 2000.

\bibitem[Babel et~al.(2010)Babel, Chuvaeva, Klin, and
  Pasechnik]{babel2010algebraic}
L.~Babel, I.~V. Chuvaeva, M.~Klin, and D.~V. Pasechnik.
\newblock Algebraic combinatorics in mathematical chemistry. {M}ethods and
  algorithms. {II}. {P}rogram implementation of the {W}eisfeiler-{L}eman
  algorithm.
\newblock \emph{arXiv preprint arXiv:1002.1921}, 2010.

\bibitem[Bachoc et~al.(2012)Bachoc, Gijswijt, Schrijver, and
  Vallentin]{bachoc2012invariant}
C.~Bachoc, D.~C. Gijswijt, A.~Schrijver, and F.~Vallentin.
\newblock Invariant semidefinite programs.
\newblock In \emph{Handbook on semidefinite, conic and polynomial
  optimization}, pages 219--269. Springer, 2012.

\bibitem[Bhatia(2009)]{bhatia2009positive}
R.~Bhatia.
\newblock \emph{Positive definite matrices}.
\newblock Princeton university press, 2009.

\bibitem[B{\"o}di et~al.(2011)B{\"o}di, Grundh{\"o}fer, and
  Herr]{bodi2011symmetries}
R.~B{\"o}di, T.~Grundh{\"o}fer, and K.~Herr.
\newblock Symmetries of linear programs.
\newblock \emph{Note di Matematica}, 30\penalty0 (1):\penalty0 129--132, 2011.

\bibitem[Borwein and Wolkowicz(1981)]{borwein1981regularizing}
J.~Borwein and H.~Wolkowicz.
\newblock Regularizing the abstract convex program.
\newblock \emph{Journal of Mathematical Analysis and Applications}, 83\penalty0
  (2):\penalty0 495--530, 1981.

\bibitem[Caluza~Machado and de~Oliveira~Filho(2018)]{caluza2018improving}
F.~Caluza~Machado and F.~M. de~Oliveira~Filho.
\newblock Improving the semidefinite programming bound for the kissing number
  by exploiting polynomial symmetry.
\newblock \emph{Experimental Mathematics}, 27\penalty0 (3):\penalty0 362--369,
  2018.

\bibitem[{de}~Klerk(2010)]{de2010exploiting}
E.~{de}~Klerk.
\newblock Exploiting special structure in semidefinite programming: A survey of
  theory and applications.
\newblock \emph{European Journal of Operational Research}, 201\penalty0
  (1):\penalty0 1--10, 2010.

\bibitem[de~Klerk and Sotirov(2009)]{de2009new}
E.~de~Klerk and R.~Sotirov.
\newblock A new library of structured semidefinite programming instances.
\newblock \emph{Optimization Methods \& Software}, 24\penalty0 (6):\penalty0
  959--971, 2009.

\bibitem[de~Klerk et~al.(2011)de~Klerk, Dobre, and Pasechnik]{de2011numerical}
E.~de~Klerk, C.~Dobre, and D.~V. Pasechnik.
\newblock Numerical block diagonalization of matrix*-algebras with application
  to semidefinite programming.
\newblock \emph{Mathematical programming}, 129\penalty0 (1):\penalty0 91--111,
  2011.

\bibitem[Dobre and Vera(2015)]{dobre2015exploiting}
C.~Dobre and J.~Vera.
\newblock Exploiting symmetry in copositive programs via semidefinite
  hierarchies.
\newblock \emph{Mathematical Programming}, 151\penalty0 (2):\penalty0 659--680,
  2015.

\bibitem[Drusvyatskiy and Wolkowicz(2017)]{drusvyatskiy2017many}
D.~Drusvyatskiy and H.~Wolkowicz.
\newblock The many faces of degeneracy in conic optimization.
\newblock \emph{arXiv preprint arXiv:1706.03705}, 2017.

\bibitem[Eberly and Giesbrecht(1996)]{eberly1996efficient}
W.~Eberly and M.~Giesbrecht.
\newblock Efficient decomposition of associative algebras.
\newblock In \emph{Proceedings of the 1996 international symposium on Symbolic
  and algebraic computation}, pages 170--178. ACM, 1996.

\bibitem[Faraut and Kor{\'a}nyi(1994)]{faraut1994analysis}
J.~Faraut and A.~Kor{\'a}nyi.
\newblock \emph{Analysis on symmetric cones}.
\newblock Oxford university press, 1994.

\bibitem[Farenick(2012)]{farenick2012algebras}
D.~Farenick.
\newblock \emph{Algebras of Linear Transformations}.
\newblock Universitext. Springer New York, 2012.
\newblock ISBN 9781461300977.

\bibitem[Fawzi and Parrilo(2014)]{fawzi2014self}
H.~Fawzi and P.~A. Parrilo.
\newblock Self-scaled bounds for atomic cone ranks: applications to nonnegative
  rank and cp-rank.
\newblock \emph{arXiv preprint arXiv:1404.3240}, 2014.

\bibitem[Faybusovich(1997)]{faybusovich1997linear}
L.~Faybusovich.
\newblock Linear systems in {J}ordan algebras and primal-dual interior-point
  algorithms.
\newblock \emph{Journal of computational and applied mathematics}, 86\penalty0
  (1):\penalty0 149--175, 1997.

\bibitem[Fujisawa et~al.(2002)Fujisawa, Kojima, Nakata, and
  Yamashita]{fujisawa2002sdpa}
K.~Fujisawa, M.~Kojima, K.~Nakata, and M.~Yamashita.
\newblock Sdpa (semidefinite programming algorithm) user’s manual—version
  6.2. 0.
\newblock \emph{Department of Mathematical and Com-puting Sciences, Tokyo
  Institute of Technology. Research Reports on Mathematical and Computing
  Sciences Series B: Operations Research}, 2002.

\bibitem[Gatermann and Parrilo(2004)]{symSoS_Gatermann200495}
K.~Gatermann and P.~A. Parrilo.
\newblock Symmetry groups, semidefinite programs, and sums of squares.
\newblock \emph{Journal of Pure and Applied Algebra}, 192\penalty0
  (1--3):\penalty0 95--128, 2004.
\newblock ISSN 0022-4049.
\newblock \doi{10.1016/j.jpaa.2003.12.011}.
\newblock URL
  \url{http://www.sciencedirect.com/science/article/pii/S0022404904000131}.

\bibitem[Gijswijt(2010)]{gijswijt2010matrix}
D.~Gijswijt.
\newblock Matrix algebras and semidefinite programming techniques for codes.
\newblock \emph{arXiv preprint arXiv:1007.0906}, 2010.

\bibitem[Grohe et~al.(2014)Grohe, Kersting, Mladenov, and
  Selman]{grohe2014dimension}
M.~Grohe, K.~Kersting, M.~Mladenov, and E.~Selman.
\newblock Dimension reduction via colour refinement.
\newblock In \emph{Algorithms-ESA 2014}, pages 505--516. Springer, 2014.

\bibitem[Hanche-Olsen and St{\o}rmer(1984)]{hanche1984jordan}
H.~Hanche-Olsen and E.~St{\o}rmer.
\newblock \emph{Jordan operator algebras}, volume~21.
\newblock Pitman Advanced Publishing Program, 1984.

\bibitem[Higman(1987)]{higman1987coherent}
D.~Higman.
\newblock Coherent algebras.
\newblock \emph{Linear Algebra and its Applications}, 93:\penalty0 209--239,
  1987.

\bibitem[Idel(2013)]{idel2013structure}
M.~Idel.
\newblock \emph{On the structure of positive maps}.
\newblock Technical University of Munich, 2013.

\bibitem[Maehara and Murota(2010)]{maehara2010numerical}
T.~Maehara and K.~Murota.
\newblock A numerical algorithm for block-diagonal decomposition of
  matrix*-algebras with general irreducible components.
\newblock \emph{Japan journal of industrial and applied mathematics},
  27\penalty0 (2):\penalty0 263--293, 2010.

\bibitem[Margot(2003)]{margot2003exploiting}
F.~Margot.
\newblock Exploiting orbits in symmetric ilp.
\newblock \emph{Mathematical Programming}, 98\penalty0 (1-3):\penalty0 3--21,
  2003.

\bibitem[Mittelmann(2003)]{mittelmann2003independent}
H.~D. Mittelmann.
\newblock An independent benchmarking of sdp and socp solvers.
\newblock \emph{Mathematical Programming}, 95\penalty0 (2):\penalty0 407--430,
  2003.

\bibitem[N{\'e}meth and N{\'e}meth(2014)]{nemeth2014lattice}
A.~N{\'e}meth and S.~N{\'e}meth.
\newblock Lattice-like subsets of {E}uclidean {J}ordan algebras.
\newblock \emph{arXiv preprint arXiv:1401.3581}, 2014.

\bibitem[Nesterov et~al.(1994)Nesterov, Nemirovskii, and
  Ye]{nesterov1994interior}
Y.~Nesterov, A.~Nemirovskii, and Y.~Ye.
\newblock \emph{Interior-point polynomial algorithms in convex programming},
  volume~13.
\newblock SIAM, 1994.

\bibitem[Packard and Doyle(1993)]{packard1993complex}
A.~Packard and J.~Doyle.
\newblock The complex structured singular value.
\newblock \emph{Automatica}, 29\penalty0 (1):\penalty0 71--109, 1993.

\bibitem[Papachristodoulou et~al.(2013)Papachristodoulou, Anderson, Valmorbida,
  Prajna, Seiler, and Parrilo]{papachristodoulou2013sostools}
A.~Papachristodoulou, J.~Anderson, G.~Valmorbida, S.~Prajna, P.~Seiler, and
  P.~Parrilo.
\newblock {SOSTOOLS} version 3.00 sum of squares optimization toolbox for
  {MATLAB}.
\newblock \emph{arXiv preprint arXiv:1310.4716}, 2013.

\bibitem[Pataki(2013)]{pataki2013simple}
G.~Pataki.
\newblock Strong duality in conic linear programming: facial reduction and
  extended duals.
\newblock \emph{Computational and Analytical Mathematics}, pages 613--634,
  2013.

\bibitem[Pataki and Schmieta(1999)]{pataki1999dimacs}
G.~Pataki and S.~Schmieta.
\newblock The {DIMACS} library of semidefinite-quadratic-linear programs.
\newblock Available at http://dimacs.rutgers.edu/Challenges/Seventh/Instances,
  1999.

\bibitem[Permenter(2018)]{permenterThesis2017}
F.~Permenter.
\newblock \emph{Reduction methods in semidefinite and conic optimization}.
\newblock PhD thesis, MIT, 2018.
\newblock URL \url{http://hdl.handle.net/1721.1/114005}.

\bibitem[Permenter and Parrilo(2015)]{permenterCoord2015}
F.~Permenter and P.~A. Parrilo.
\newblock Finding sparse, equivalent {SDPs} via minimal-coordinate-projections.
\newblock In \emph{IEEE 54th Annual Conference on Decision and Control (CDC)}.
  IEEE, 2015.

\bibitem[Schrijver(1979)]{schrijver1979comparison}
A.~Schrijver.
\newblock A comparison of the {D}elsarte and {L}ov{\'a}sz bounds.
\newblock \emph{Information Theory, IEEE Transactions on}, 25\penalty0
  (4):\penalty0 425--429, 1979.

\bibitem[Seiler(2013)]{seiler2013sosopt}
P.~Seiler.
\newblock {SOSOPT}: A toolbox for polynomial optimization.
\newblock \emph{arXiv preprint arXiv:1308.1889}, 2013.

\bibitem[St{\o}rmer(2013)]{stormer2012positive}
E.~St{\o}rmer.
\newblock \emph{Positive linear maps of operator algebras}.
\newblock Springer Science \& Business Media, 2013.

\bibitem[St{\o}rmer and Effros(1979)]{stormer1979positive}
E.~St{\o}rmer and E.~G. Effros.
\newblock Positive projections and {J}ordan structure in operator algebras.
\newblock \emph{Mathematica Scandinavica}, 45:\penalty0 127--138, 1979.

\bibitem[Sturm(1999)]{sturm1999using}
J.~F. Sturm.
\newblock Using {SeDuMi} 1.02, a {MATLAB} toolbox for optimization over
  symmetric cones.
\newblock \emph{Optimization methods and software}, 11\penalty0 (1-4):\penalty0
  625--653, 1999.

\bibitem[Vallentin(2009)]{vallentin2009symmetry}
F.~Vallentin.
\newblock Symmetry in semidefinite programs.
\newblock \emph{Linear Algebra and Its Applications}, 430\penalty0
  (1):\penalty0 360--369, 2009.

\bibitem[Wedderburn(1908)]{wedderburn1908hypercomplex}
J.~M. Wedderburn.
\newblock On hypercomplex numbers.
\newblock \emph{Proceedings of the London Mathematical Society}, 2\penalty0
  (1):\penalty0 77--118, 1908.

\bibitem[Weisfeiler(1977)]{weisfeiler1977construction}
B.~Weisfeiler.
\newblock \emph{On construction and identification of graphs}.
\newblock Springer, 1977.

\end{thebibliography}

}

\section{Appendix}

\subsection{Proof of Theorem~\ref{thm:posjor}}
We now prove  Theorem~\ref{thm:posjor}, which stated that a subspace $\mathcal{S} \subseteq \mathbb{S}^n$
is a Jordan subalgebra if and only if its orthogonal projection $P_{\mathcal{S}}$ is unital and positive.  
Analogues for complex Jordan algebras are well known; see~\cite{stormer2012positive}~\cite{stormer1979positive}  and also the thesis~\cite{idel2013structure}.  One direction is also shown in~\cite{nemeth2014lattice}.    The converse direction is shown in part by translating an argument of~\cite{stormer2012positive} from the complex to real case.    
Since they are short and self-contained, we  give full proofs of both directions. 

To begin, we need the following lemma relating invariance under squaring  to eigenvalue decompositions. 
\begin{lem}\label{lem:squarechar}
	For a non-zero $X \in \mathbb{S}^n$,  let $E_X \subset \mathbb{S}^n$ be the set of pairwise orthogonal idempotent matrices  for which  
	\[
	X =  \sum_{ E \in E_X} \lambda_E E,
	\]
        where the range of $E \in E_X$ is an eigenspace of $X$ and ${\{\lambda_E\}}_{E \in E_X}$ is the set of non-zero (distinct) eigenvalues of $X$.  For a subspace $\mathcal{S}  \subseteq \mathbb{S}^n$, the following are equivalent.
	\begin{enumerate}
		\item $\mathcal{S}$ contains the set  $E_X$ for all non-zero $X \in \mathcal{S}$.
		\item $\mathcal{S}$ is invariant under   squaring, i.e., $\mathcal{S} \supseteq \{ X^2 : X \in \mathcal{S} \}$.
	\end{enumerate}
	\begin{proof}
                That statement one implies two is immediate given that  $X^2 =
                \sum_{ E \in E_X } \lambda^2_E E$. Conversely,  suppose $X$ has
                non-zero eigenvalue $\lambda$  of maximum magnitude. Then, if
                statement two holds, the idempotent $\hat E = \lim_{n
                \rightarrow \infty}{(|\lambda|^{-1}  X  )}^{2n}$ is contained in
                $\mathcal{S}$ and has range equal to an eigenspace or, if
                $\pm|\lambda|$ are both eigenvalues, the sum of two
                eigenspaces. Replacing $X$ with $X-\lambda\hat E$ and iterating
                yields a set of idempotents whose span contains $E_X$;
                moreover, this set is contained in $\mathcal{S}$.
	\end{proof}
\end{lem}
\noindent  We now use this lemma and the mentioned argument of~\cite{stormer2012positive} to prove Theorem~\ref{thm:posjor}

		To prove ($2 \Rightarrow 1)$,  consider  $X \succeq 0$ and suppose $P_{\mathcal{S}}(X)$ is non-zero. For a non-zero eigenvalue $\lambda_E$ of $P_{\mathcal{S}}(X)$, let $E \in \mathbb{S}^n$ denote the idempotent with range equal to the associated eigenspace.  If (2) holds, then Lemma~\ref{lem:squarechar} implies  $P_{\mathcal{S}}(E) = E$. 
		Hence,
		\[
		0  \le  \ip{E}{X}  = \ip{P_{\mathcal{S}}(E)}{X}  =  \ip{E}{P_{\mathcal{S}}(X)}  =  \lambda_E \norm{E}^2.
		\]
		We conclude the eigenvalues of $P_{\mathcal{S}}(X)$ are non-negative, i.e., that $P_{\mathcal{S}}(X) \succeq 0$.
		To show the unitality condition, let $Z$ be a matrix in $\mathcal{S}$ of maximum rank and let
		\[
		\hat E= \sum_{ E \in E_Z} E.
		\]
		For all $X \in \mathcal{S}$, it holds that $t \hat E \succeq X^2$ for some $t > 0$. This shows the range of $\hat E$ contains the range of $X^2$ and hence the range of $X$.  It follows $\hat E X=X$.

                To prove ($1 \Rightarrow 2)$, suppose the unit element $E$ has rank $r$. Then we can find an orthogonal matrix  $Q = (Q_1, Q_2) \in \mathbb{R}^{n \times n}$ for which $E=Q_{1} Q_{1}^T$ and
		\[
		\mathcal{S} = \left\{ Q \left( \begin{array}{cc} X & 0 \\ 0 & 0 \end{array} \right) Q^T : X \in \mathcal{\hat S} \subseteq \mathbb{S}^{r}\right\},
		\]
                where $\mathcal{\hat S} := Q_{1}^T \mathcal{S} Q_{1}$. Further, the projection $P_{\mathcal{S}}$ satisfies
		\[
                  P_{\mathcal{S}}(X) = Q_{1} Q_{1}^T P_{\mathcal{\hat S}  }(X) Q_{1} Q_{1}^T
		\] 
		where $P_{\mathcal{\hat S}  } :  \mathbb{S}^{r} \rightarrow \mathbb{S}^{r}$ is the orthogonal projection onto $\mathcal{\hat S}$.
                It follows that if $\mathcal{\hat S}$ is invariant under   squaring, so is $\mathcal{S}$, and if $P_{\mathcal{S}}$ is positive, so is $P_{\mathcal{\hat S}  }$.  
		Hence, Statement 2 follows by showing $\mathcal{\hat S}$ is invariant under   squaring.
                
                We show this  applying the argument from~\cite[Theorem 2.2.2]{stormer2012positive} and using the fact $\mathcal{\hat S}$ contains the identity matrix of order $r$. Dropping the subscript $\mathcal{\hat S}$ from $P_{\mathcal{\hat S}}$, we first note since  $P$ is positive and $P(I)=I$, it satisfies the Kadison inequality, which states  $P(X^2) -P(X)P(X)  \succeq 0$ for all $X \in \mathbb{S}^r$ (e.g., Theorem 2.3.4 of~\cite{bhatia2009positive}).  Hence, for $X$ in the range of $P$
		\[
		P(X^2) -X^2 \succeq 0. 
		\]
		Letting $Z=P(X^2) -X^2$ and taking the trace shows
		\[
		\trace Z = \ip{I}{Z} =  \ip{P(I)}{Z} = \ip{I}{P(Z)} = \trace \big(P^2(X^2) -P(X^2)\big) = \trace \big( P(X^2) -P(X^2) \big) = 0.
		\]
		Since $Z \succeq 0$, we conclude $Z=0$, i.e., that $P(X^2)=X^2$. Therefore $X^2$ is in the range of $P$.

\subsection{Invariant affine sets of projections}\label{sec:invarAffine}
Recall Condition~\ref{cond:CI}-\ref{itm:primaff} and Condition~\ref{cond:CI}-\ref{itm:dualaff} require  invariance of the affine sets $\primaff$ and $\dualaff$ under the projection $P_{\mathcal{S}}$. We now
prove the characterization of these conditions provided by Lemma~\ref{lem:subspacechar}.
\begin{lem} 
	For an affine set  $Y + \mathcal{L}$,  let 
	$Y_{\mathcal{L}^{\perp}} \in \mathbb{S}^n$  denote the projection  of $Y\in \mathbb{S}^n$   onto the subspace    $\mathcal{L}^{\perp}$.  Let $P_{\mathcal{S}} : \mathbb{S}^n \rightarrow \mathbb{S}^n$ denote the
	orthogonal projection onto a subspace $\mathcal{S}$ of $\mathbb{S}^n$. The following statements are equivalent.
	\begin{enumerate}
		\item $P_{\mathcal{S}}(Y  + \mathcal{L}) \subseteq Y  + \mathcal{L}$  
		\item $P_{\mathcal{S}}(Y_{\mathcal{L}^{\perp}}) = Y_{\mathcal{L}^{\perp}}$ and $P_{\mathcal{S}}(\mathcal{L}) \subseteq \mathcal{L}$
	\end{enumerate}
	\begin{proof}
		 To begin, first note $P_{\mathcal{S}}$---being an orthogonal projection---is a contraction with respect to the   Frobenius norm $\| X \|_{F}$ (recalling our use of the trace inner-product); further, $Y_{\mathcal{L}^{\perp}}$ is the unique minimizer of this norm over $Y  + \mathcal{L}$.  Hence, if   $P_{\mathcal{S}}(Y  + \mathcal{L}) \subseteq Y  + \mathcal{L}$, then $P_{\mathcal{S}}(Y_{\mathcal{L}^{\perp}}) = Y_{\mathcal{L}^{\perp}}$; in addition, since $Y  + \mathcal{L} = Y_{\mathcal{L}^{\perp}}  + \mathcal{L}$,
		\[
		Y_{\mathcal{L}^{\perp}}  + P_{\mathcal{S}}(\mathcal{L}) =   P_{\mathcal{S}}(  Y_{\mathcal{L}^{\perp}} +\mathcal{L})  \subseteq Y_{\mathcal{L}^{\perp}}  + \mathcal{L},
		\]
		which implies $P_{\mathcal{S}}(\mathcal{L}) \subseteq \mathcal{L}$. The converse direction is obvious given that  $Y  + \mathcal{L} = Y_{\mathcal{L}^{\perp}}  + \mathcal{L}$.
		
	\end{proof}
\end{lem}
\noindent If we apply the previous lemma to both the primal and dual affine sets we obtain the conditions $P_{\mathcal{S}}(\mathcal{L}) \subseteq \mathcal{L}$ and $P_{\mathcal{S}}(\mathcal{L}^{\perp}) \subseteq \mathcal{L}^{\perp}$. However, Lemma~\ref{lem:subspacechar} only contains one of these conditions, since they turn out to be equivalent. Consider the following.
\begin{lem}{~\cite[Proposition 3.8]{farenick2012algebras} }
	Let $P_{\mathcal{L}} : \mathbb{S}^n \rightarrow \mathbb{S}^n$ and $P_{\mathcal{S}}: \mathbb{S}^n \rightarrow \mathbb{S}^n$ denote the orthogonal
	projections onto subspaces $\mathcal{L}$ and $\mathcal{S}$ of $\mathbb{S}^n$.  The following four statements are equivalent.
	\begin{multicols}{2}
		\begin{itemize}
			\item $\mathcal{L}$ is an invariant subspace of $P_{\mathcal{S}}$
			\item $\mathcal{L}^{\perp}$ is an invariant subspace of $P_{\mathcal{S}}$
			\item $\mathcal{S}$ is an invariant subspace of $P_{\mathcal{L}}$
			\item $\mathcal{S}^{\perp}$ is an invariant subspace of $P_{\mathcal{L}}$
		\end{itemize}
	\end{multicols}

\end{lem} 
\noindent  Combining  these   two lemmas proves Lemma~\ref{lem:subspacechar}.

\subsection{Linear images of self-dual cones}
\renewcommand{\coneName}{\mathcal{K}}
The following was used to prove Proposition~\ref{prop:equiv2}.  
\begin{lem}\label{lem:sdimage}
	Let $\mathcal{W}$ and $\mathcal{V}$ be   inner-product spaces and $\coneNameTwo \subseteq \mathcal{V} $ and $\coneName  \subseteq \mathcal{W}$  self-dual convex cones. Let $T : \mathcal{V} \rightarrow \mathcal{W}$ be a injective linear map with adjoint $T^* : \mathcal{W} \rightarrow \mathcal{V}$.  If  $\coneName = T(\coneNameTwo)$, then
	$T^*T(\coneNameTwo) = \coneNameTwo$.
	
	\begin{proof}
		
		For all $x,y \in \coneNameTwo$, 
		\[
		\langle  T^*T(x),  y \rangle = \langle  T(x),  T(y) \rangle \ge 0
		\]
                by self-duality of $\coneName$. By self-duality of
                $\coneNameTwo$, we conclude $T^*T(x) \in \coneNameTwo$. On the
                other hand, since $T^*$ is surjective, we have for any $x \in
                \coneNameTwo$ existence of $w\in \mathcal{V}$ for which $x=T^*w$. Further, for
                all $y \in \coneNameTwo$,
		\[
		0 \le \langle T^*w, y\rangle =  \langle w, T y\rangle
		\]
		which, since $\coneName = T(\coneNameTwo)$, shows $w \in \coneName$. Hence, $w=Tz$ for $z \in \coneNameTwo$, showing $x = T^*Tz$.
	\end{proof}
	
\end{lem}

\end{document}